\documentclass[11pt]{amsart}

\title{Cramér distance and discretizations of circle expanding maps II: simulations}
\date{\today}

\thanks{This project was partially supported by a PEPS/CNRS project and the ANR CODYS}

\author{Pierre-Antoine Guih\'eneuf}
\address{Pierre-Antoine Guih\'eneuf: Sorbonne Université and Université de Paris, CNRS, IMJ-PRG, F-75005 Paris, France.}
\email{pierre-antoine.guiheneuf@imj-prg.fr}

\author{Maurizio Monge}
\address{Maurizio Monge: Instituto de Matem{\`a}tica da UFRJ, Av. Athos da Silveira Ramos 149, Centro de Tecnologia, Bloco C Cidade Univesit{\`a}ria, Ilha do Fund{\~a}o, Caixa Postal 68530 21941-909, Rio de Janeiro, RJ, Brasil}
\email{maurizio.monge@im.ufrj.br}

\usepackage[english]{babel}
\usepackage{amsfonts}
\usepackage{amsmath}
\usepackage{amssymb}
\usepackage{stmaryrd}
\usepackage{amsthm}
\usepackage{enumitem}
\usepackage{pgf,tikz}
\usetikzlibrary{decorations.pathreplacing, shapes.multipart, arrows, matrix}
\usepackage[margin=0pt]{caption}
\usepackage[pdftex,colorlinks=true,linkcolor=blue,citecolor=blue,urlcolor=blue]{hyperref}

\newtheorem{lemme}{Lemma}
\newtheorem{lemma}[lemme]{Lemma}

\newtheorem{theoreme}[lemme]{Theorem}
\newtheorem{prop}[lemme]{Proposition}
\newtheorem{coro}[lemme]{Corollary}
\newtheorem{conj}[lemme]{Conjecture}

\newtheorem{ques}[lemme]{Question}

\theoremstyle{definition}

\newtheorem*{moral}{\textbf{Moral}}

\theoremstyle{remark}
\newtheorem{rem}[lemme]{Remark}

\newcommand{\C}{\mathbf{C}}

\newcommand{\N}{\mathbf{N}}

\newcommand{\R}{\mathbf{R}}

\newcommand{\Z}{\mathbf{Z}}
\newcommand{\D}{\mathcal{D}}

\newcommand{\Leb}{\operatorname{Leb}}

\newcommand{\ud}{\,\mathrm{d}}

\newcommand{\card}{\operatorname{Card}}

\newcommand{\Sp}{\mathbf{S}}

\newcommand{\Disc}{\operatorname{d_C}}

\newcommand\numberthis{\addtocounter{equation}{1}\tag{\theequation}}

\allowdisplaybreaks

\begin{document}

\begin{abstract}
This paper presents some numerical experiments in relation with the theoretical study of the ergodic short-term behaviour of discretizations of expanding maps done in \cite{paper1}.

Our aim is to identify the phenomena driving the evolution of the distance between the $t$-th iterate of Lebesgue measure by the dynamics $f$ and the $t$-th iterate of the uniform measure on the grid of order $N$ by the discretization on this grid. Based on numerical simulations we propose some conjectures on the effects of numerical truncation from the ergodic viewpoint.
\end{abstract}

\maketitle

\setcounter{tocdepth}{1}
\tableofcontents

\section{Introduction}

\begin{figure}[ht]
\begin{minipage}{0.48\linewidth}   
\centering \includegraphics[width=\linewidth]{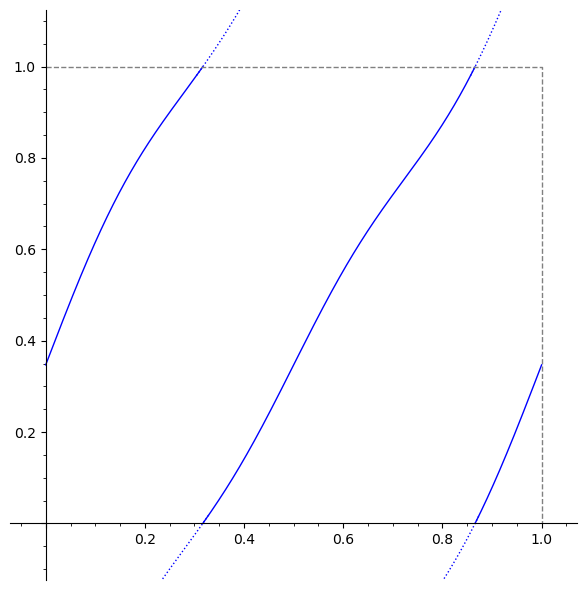}
\caption{\label{FigMap} Graph of the studied expanding map (see  Subsection~\ref{SecMaps}).}
\end{minipage}
\hfill
\begin{minipage}{.48\linewidth}
\centering \includegraphics[width=\linewidth]{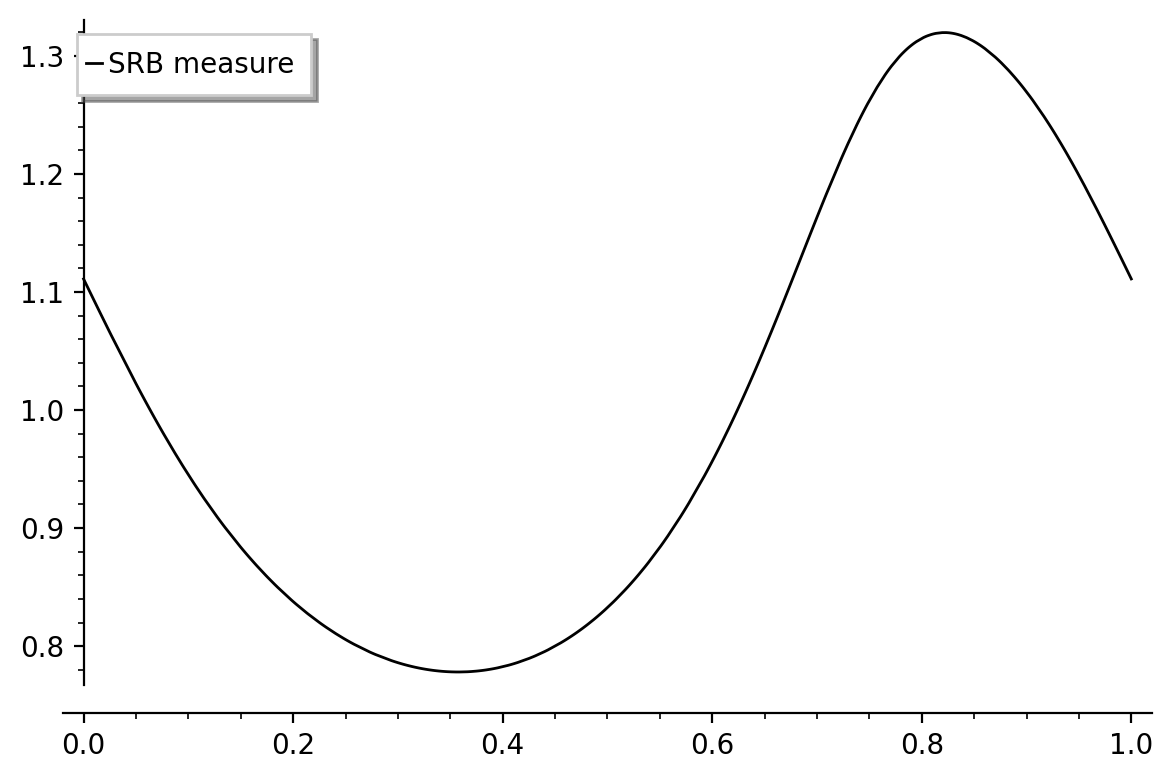}
\caption{\label{FigDensSRB} Density of the SRB measure associated to the map of Figure \ref{FigMap}.}
\end{minipage}
\end{figure}

This article is the experimental part of a series of two papers aiming to understand the ergodic behaviour of discretizations of circle expanding maps (see \cite{paper1}). By expanding map of the circle $\Sp^1 = \R/\Z$ we mean a $C^{r}$ map $f : \Sp^1\to \Sp^1$ ($r>1$) such that $f'(x)>1$ for any $x\in\Sp^1$ (see Figure \ref{FigMap} for an example of such a map, described in Subsection~\ref{SecMaps}). Note that these assumptions force the map to be of degree $d\ge 2$.

We identify the circle $\Sp^1 \simeq \R/\Z$ with its fundamental domain $[0,1]$, and endow it with discretization grids, of parameter $N> 0$
\[E_N = \left\{\frac{i}{N} \mid 0\le i <N\right\},\]
and discretization projections $P_N : \Sp^1 \to E_N$ defined by 
\[P_N(x) = \frac{i}{N} \iff x\in \left[\frac{i-\frac12}{N} , \, \frac{i+\frac12}{N} \right).\]
This allow to define the \emph{discretizations} $f_N : P_N \to P_N$ of the map $f$ by $f_N = P_N \circ f|_{E_N}$. In other words, $f_N(x)$ is obtained from $f(x)$ by projecting on the closest point of the grid $E_N$. Of course, this models what happens when the computer iterates a map using a fixed number of digits --- when $N = 2^k$, the set $E_N$ represents the set of points with at most $k$ binary places. We also set $\Leb_N$ the uniform probability measure on $E_N$.
\medskip

The basic example of expanding map $f : x\mapsto 2x$ shows that in some cases the  discretizations dynamics does not reflect the chaotic properties of the map: if $N = 2^k$, then $f_N = f|_{E_N}$ and $f_N^k(x) = 0$ for any $x\in E_N$. In other words, any point of the grid is mapped after a small number of iterations on the fixed point 0: the dynamics of $f_N$ is completely trivial.

To avoid these phenomena of resonance between the dynamics and the grid --- that one can expect to be exceptional --- one can consider \emph{generic} dynamics. A property on expanding maps is said to be generic if it is satisfied on at least a countable intersection of open and dense subsets of the space of $C^{r}$ expanding maps (for $C^r$ topology). Baire's theorem ensures that a generic property is satisfied on a dense set of dynamics.
\medskip

While some theoretical results are known about the \emph{local} dynamics of discretizations of $C^1$ generic dynamics (e.g. \cite{Gui15c}), to our knowledge, the only known result about their \emph{global} dynamics deals with the \emph{degree of recurrence} (see \cite{vladimirov2015quantized} and \cite{MR3919917}). Besides this local/global dichotomy, one can classify the discretizations' dynamics into combinatorial and ergodic properties. Whereas combinatorial properties have been the subject of numerous numerical explorations, ergodic properties have been only little studied.

In this work (together with \cite{paper1}), we intend to study the global ergodic behaviour of generic circle expanding maps discretizations.


The smoothness assumption on $f$ ensures the existence of a unique absolutly continuous invariant measure, called SRB (for Sinai-Ruelle-Bowen), which is moreover ergodic, mixing and has the property that (and this is crucial here) the measures $f_*^k(\Leb)$ converge exponentially fast to it\footnote{Meaning that the densities of these measures converge exponentially fast towards the density of $SRB$ in the $C^{r-1}$ topology}. This measure is of great importance for the ergodic study of expanding maps, and its counterpart for higher dimensional hyperbolic maps opened the way to a whole branch of the ergodic theory.
See Figure \ref{FigDensSRB} for the graph of this measure's density in the case of the map of Figure \ref{FigMap}.
\medskip

A large part of this paper will be devoted to the numerical comparison, for some expanding map $f$, of the actions of $f$ and of its discretizations $f_N$ on uniform measures. On the one hand, as said before, the iterates $f^k_*(\Leb)$ converge exponentially fast towards $SRB$. On the other hand, what happens to the measures $(f_N^k)_*(\Leb_N)$, where $\Leb_N$ denotes the uniform measure on $E_N$, is much more unclear. 

\begin{figure}[ht]
\includegraphics[width=.48\linewidth]{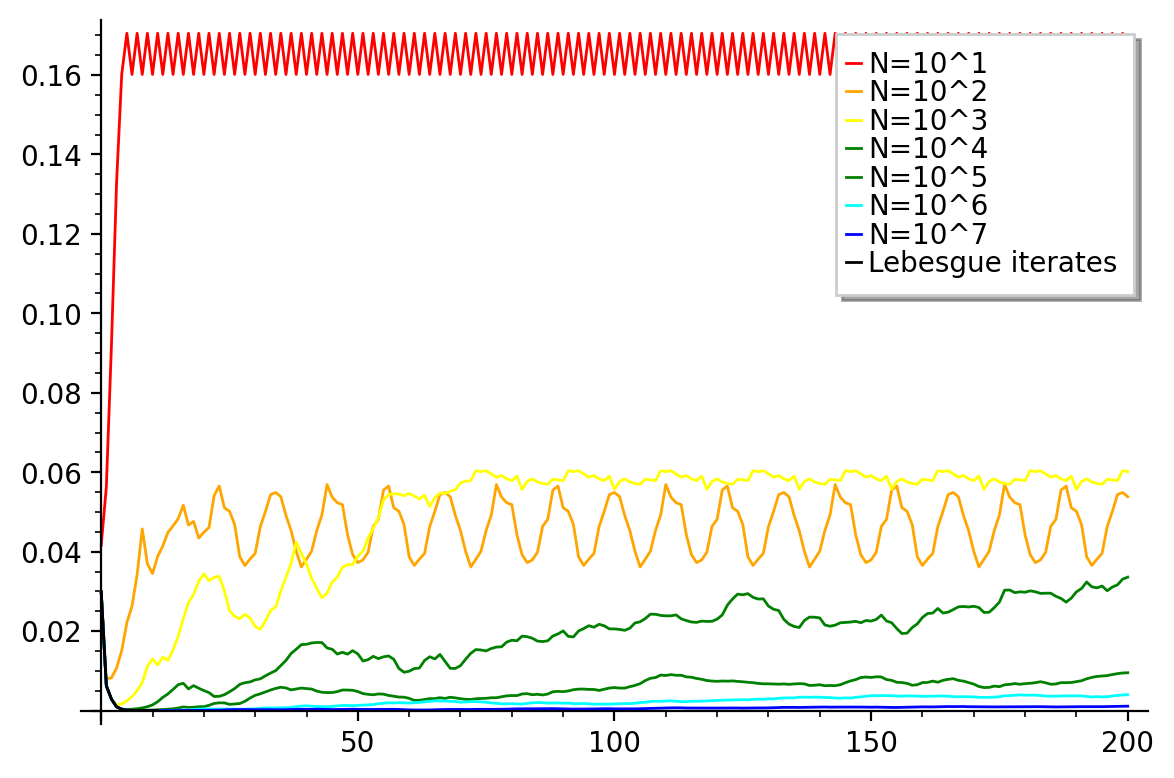}
\hfill
\centering \includegraphics[width=.48\linewidth]{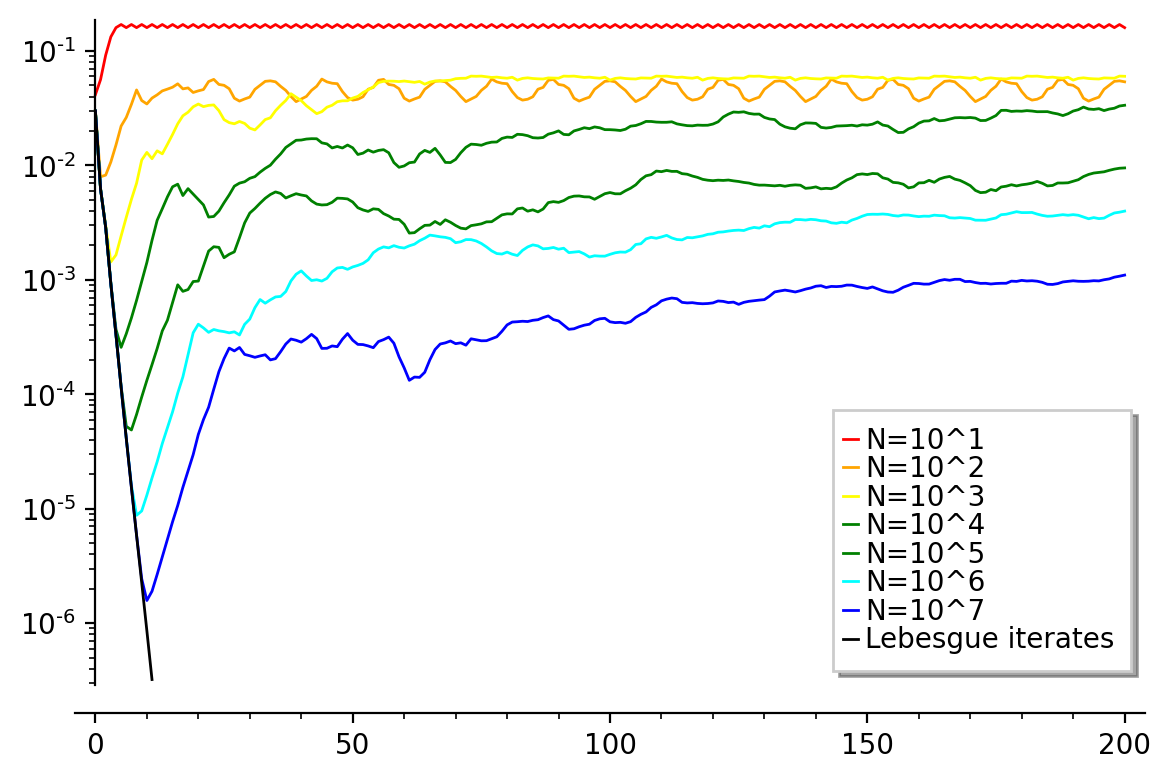}
\caption{\label{FirstSimul} Graphs of $k\mapsto\Disc(f_*^k(\Leb),SRB)$ (black) and $k\mapsto\Disc((f_N)_*^k(\Leb_N),SRB)$ for various grid sizes $N$, uniform (left) and logarithmic (right) scales.}
\end{figure}

Figure \ref{FirstSimul} shows the evolution of $\Disc(f_*^k(\Leb),SRB)$ and $\Disc((f_N)_*^k(\Leb_N),SRB)$ with $k$, for different discretization orders $N$. Here, $\Disc$ is a distance on the set of probability measures, that we call \emph{Cramér distance}, defined in Equation \eqref{DefDisc} page \pageref{DefDisc}, and which spans the weak-* topology. This distance il also called Cramér-von Mises distance, or ``the $L_2$-metrics between distribution functions'' \cite{MR1105086, zbMATH05216873}

\begin{figure}[ht]
\centering
\includegraphics[width=.48\linewidth]{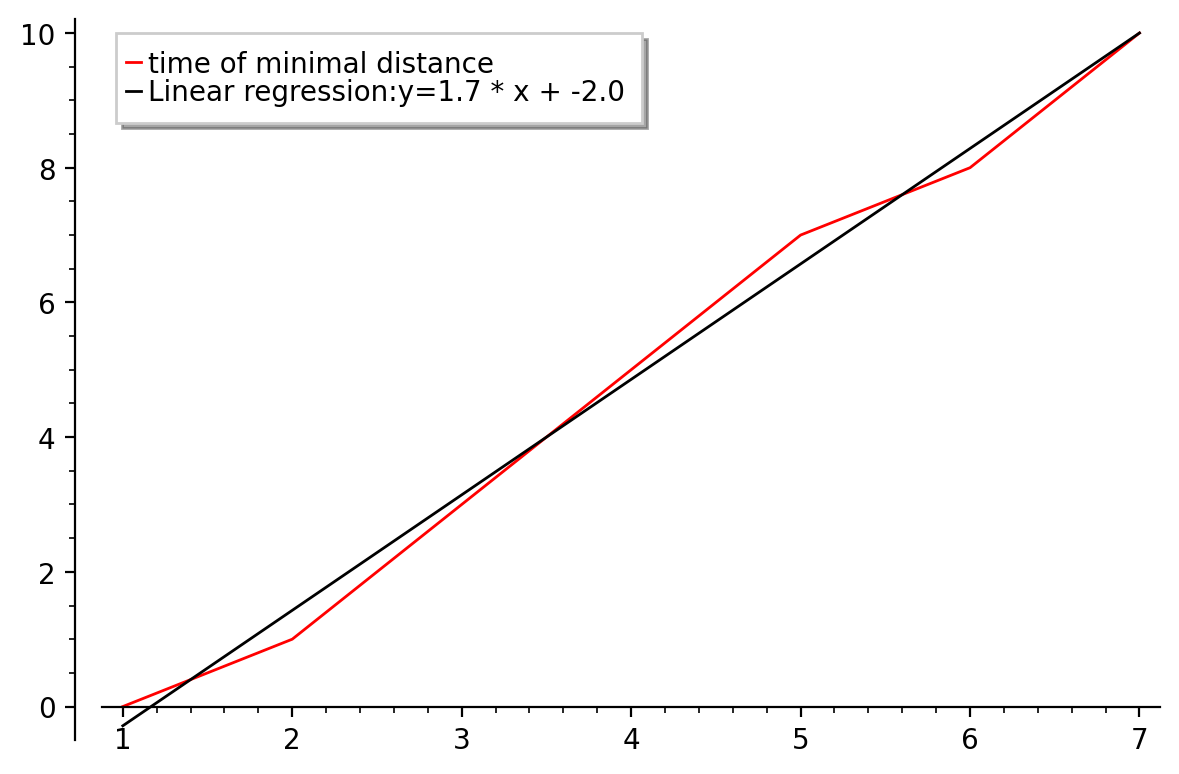}
\caption{\label{FirstSimul2} Argument of the minimum of the distance  $k\mapsto\Disc((f_N)_*^k(\Leb_N),SRB)$ (i.e. time $t_N$ for which this distance is minimal) depending on $\log_{10}(N)$, and linear regression of these values. This strongly suggests that this time $t_N$ is of the order of $\log N$.}
\end{figure}

\begin{figure}
\includegraphics[width=.48\linewidth]{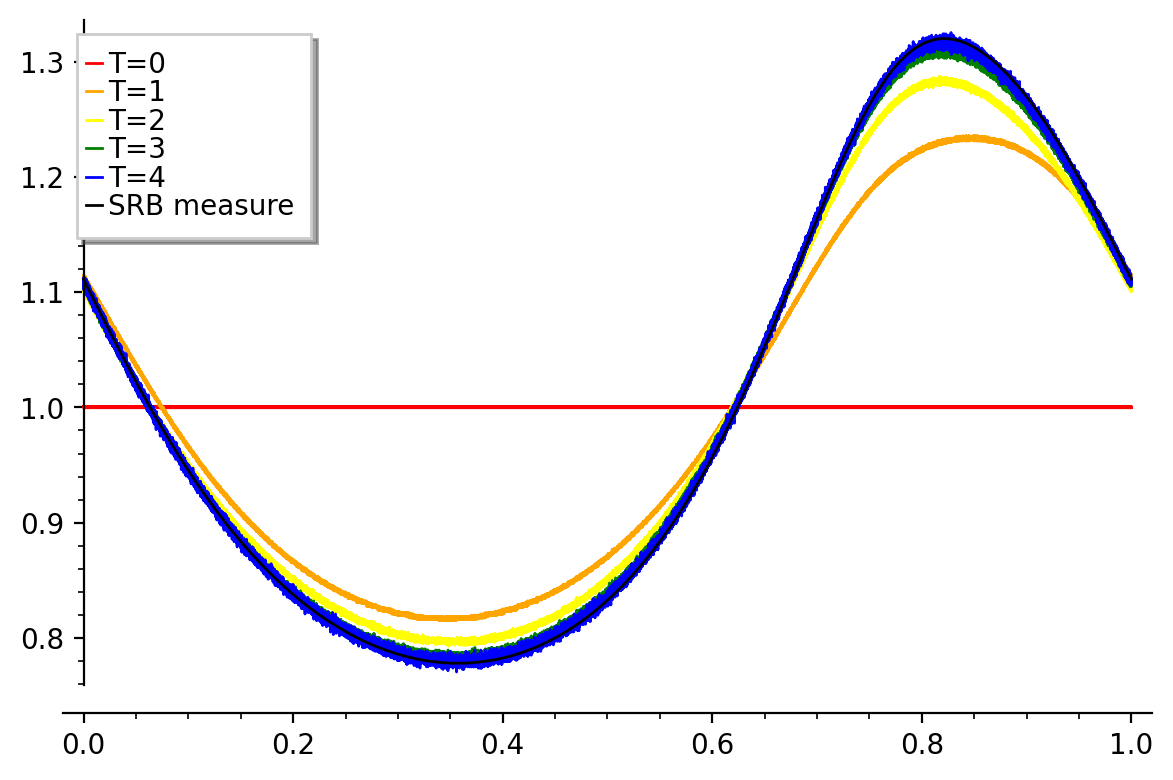}
\hfill
\includegraphics[width=.48\linewidth]{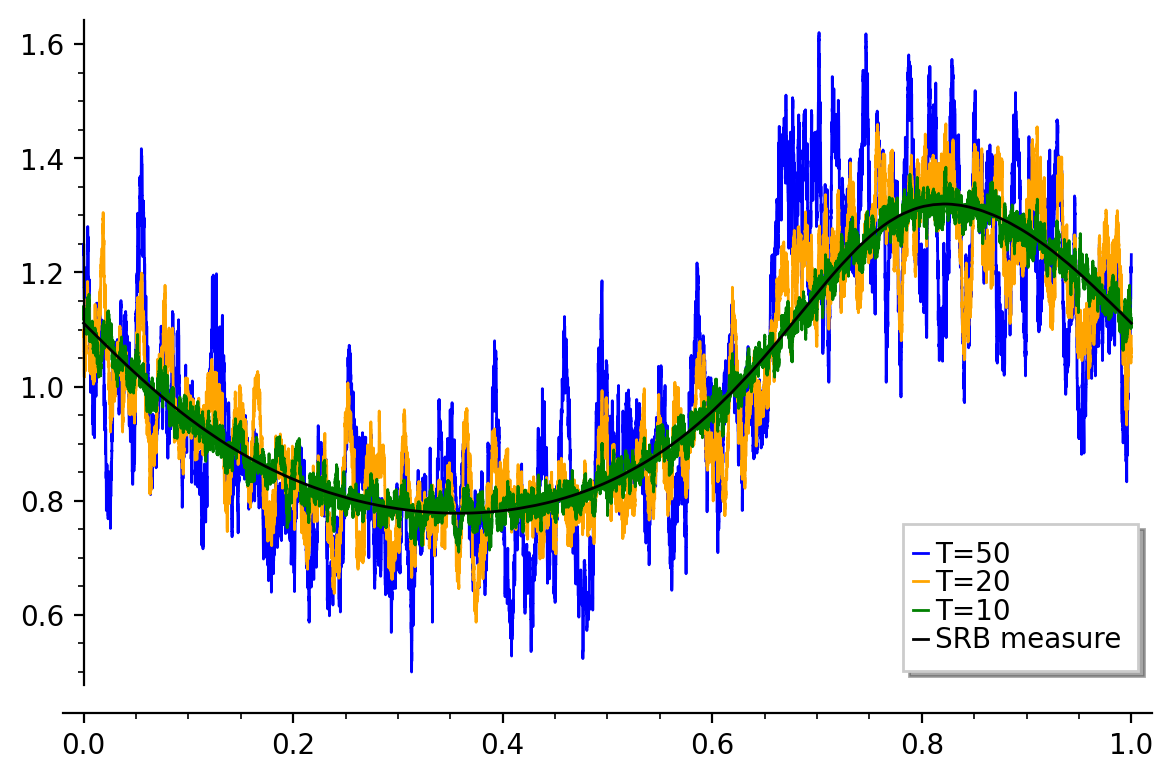}
\caption{Densities of the measures $f_N^T(\Leb_N)$ for $N=10^5$ and for different iteration times $T$, together with the density of the SRB measure. Note that as the measures $f_N^T(\Leb_N)$ are punctual, these measures are smoothed on intervals of size $500/N$ (i.e. they are convolved with the indication function of these intervals). \label{FigDensityIterates}}
\end{figure}

All the curves for the discretizations have more or less the same shape:
\begin{itemize}
\item They first decrease, up to a certain point, following quite well the corresponding curve for the actual dynamics $f$ (in black), which decreases exponentially (see also Figure~\ref{FigDensityIterates}, left, for the plot of the densities of these measures). Figure~\ref{FirstSimul2} suggests a more or less linear relation between this time of minimum of the distance and $\log N$.
\item Then they move away from the black curve and start to increase (see also Figure~\ref{FigDensityIterates}, right).
\item From a certain point, they seem to have a periodic behaviour (at least the ones for small values of the order $N$).
\end{itemize}

Let us explain the behaviour for small times. As a consequence of Theorem~\ref{MainTheo}, for any $k\ge 0$, one has:
\[(f_N)^k_*(\Leb_N) \underset{N\to +\infty}{\longrightarrow} f^k_*(\Leb).\]
This is illustrated by Figure~\ref{FigDensityIterates}. See also \cite[Theorem 12.17]{Guih-These} for a proof with effective bounds on convergence speed. Roughly speaking, the operators $(f_N)_*$, acting on invariant measures, converge towards $f_*$. Hence, the behaviour of
\[k\mapsto \Disc\big((f_N)_*^k(\Leb_N),SRB\big)\]
is the ``combination'' of the behaviours of 
\begin{align}
k & \mapsto\Disc\big((f_N)_*^k(\Leb_N),f_*^k(\Leb)\big)
\qquad \text{and}\label{EqStudiedQuant}\\
k & \mapsto\Disc\big(f_*^k(\Leb),SRB\big).\nonumber
\end{align}
The second one is well understood, as $\Disc(f_*^k(\Leb),SRB)$ tends to 0 exponentially fast in $k$ (combine Lemma~\ref{lemDiscrepDensit} with the fact that the densities converge exponentially to the density of $SRB$ in the $C^{r-1}$ topology), so we are reduced to study the first map \eqref{EqStudiedQuant}. 

\begin{figure}[ht]
\includegraphics[width=.49\linewidth]{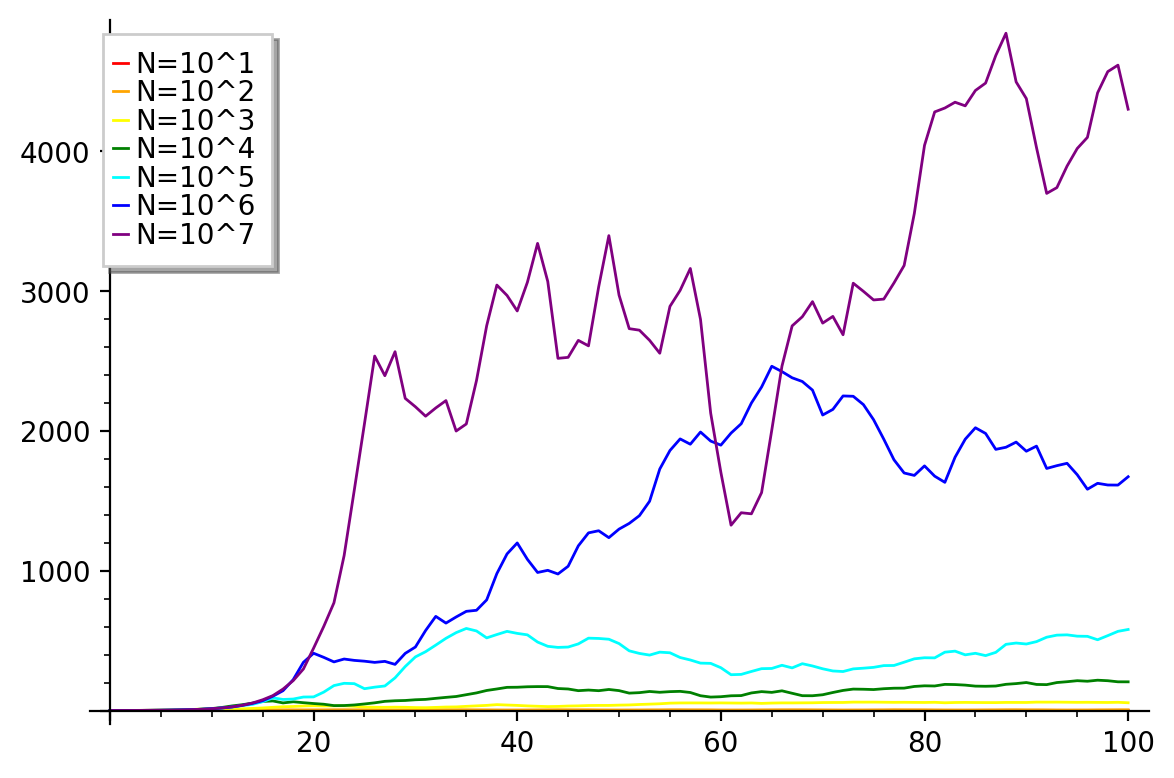}
\hfill
\centering \includegraphics[width=.49\linewidth]{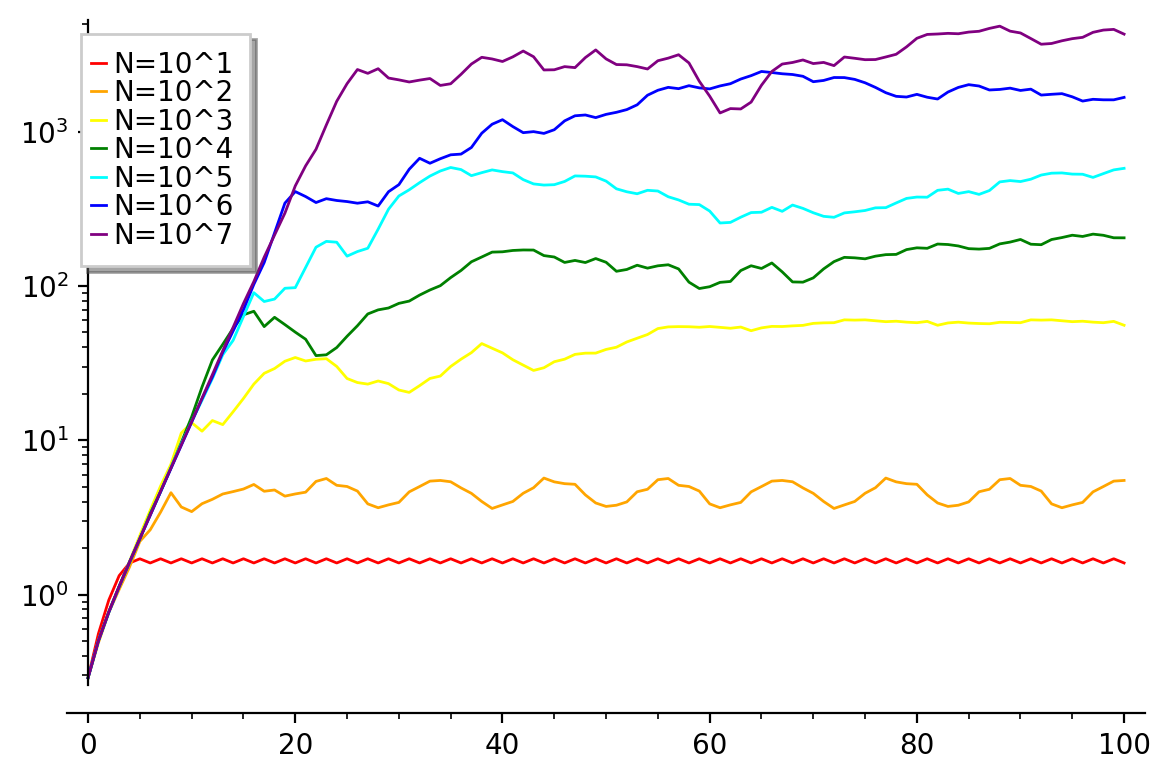}
\caption{\label{SecondSimul} Graphs of $k\mapsto N\Disc((f_N)_*^k(\Leb_N),f_*^k(\Leb))$ for various grid sizes $N$, uniform (left) and logarithmic (right) scales. The factor $N$ in front of the distance $\Disc$ is added so that the curves all start at $1/\sqrt{12}$ for the time $k=0$ (the distance $\Disc$ between Lebesgue measure and the uniform measure on $E_N$ is $1/(N\sqrt{12})$).}
\end{figure}

Figure~\ref{SecondSimul} shows simulations of the first map \eqref{EqStudiedQuant}.
On these graphics there are three distinct time regimes:
\begin{enumerate}[label=($R_\arabic*$), ref=($R_\arabic*$)]
\item\label{R1} The short-term behaviour, where the curves for discretizations seem to have a uniform behaviour: for a fixed time $k$ and the grid parameter $N$ going to infinity, they seem to converge to a curve that depends more or less exponentially on $k$ (it is close to a line on the right graph which is in logarithmic scale). For this regime we have a theoretical prediction given by the asymptotics \eqref{EqIntMieuxIntro} of Theorem~\ref{MainTheo}. We will confront this prediction with the actual simulations in the sequel.
\item\label{R2} The medium term behaviour, where the curve globally grows slowly.
\item\label{R3} The long-term or asymptotic behaviour, where the curve is periodic (at least for small values of $N$).
\end{enumerate}

These three different regimes will guide our study of \eqref{EqStudiedQuant}. More precisely, we will study these three regime separately and one after the other.
\medskip

As proved independently in \cite{MR1661120} and \cite{Flocker}, for a generic map of $\D^r(\Sp^1)$, the roundoff errors equidistribute: for a fixed time $k$, and the order $N$ going to infinity, the sequence of roundoff errors in time $k$ equidistributes in $[-1/(2N),1/(2N)]^k$. See also \cite[Proposition 3.2]{paper1} for a more precise statement, which is obtained as a byproduct of the proofs.

So at first sight, one could expect the discretizations to behave very similarly to random perturbations. More precisely, the discretization $f_N$'s global dynamics may be thought as the typical global dynamics of the random map $f_t$, acting on $N$-tuples of points of $\Sp^1$ such that each point $x$ of this tuple is randomly drawn uniformly in $[f(x)-1/(2N),f(x)+1/(2N)]$.

In fact, things are a bit more subtle, and one quickly realizes that the fact that orbits of $f_N$ can merge (i.e. that there exists distinct grid points that are eventually mapped to the same point under $f_N$) --- and hence will stay together forever --- must be an important parameter influencing the evolution of $(f_N)_*^k(\Leb_N)$. With this in mind, one can isolate (at least) four phenomena that make the action of discretizations different from that of a random map.
\begin{enumerate}[label=($P_\arabic*$), ref=($P_\arabic*$)]
\item The iterates of points always belong to $E_N$.
\item Two points of $E_N$ having the same image by $f_N$ will have identical positive orbits.
\item The local shape around $y\in\Sp^1$ of the image $f_N(E_N)$ is very similar to the one of a linearization of $f$ around the points $f^{-1}(y)$, which is a model set (see \cite{paper1}).
\item Any point eventually falls in a periodic cycle.
\end{enumerate}

Part of the paper will be devoted to the understanding of the relative effects of these phenomena on the action of discretizations on measures.

\subsection*{Some bibliographical remarks}

There are numerous numerical studies of the spatial discretization effects, but strangely only few works about these effects on ergodic properties.

A lot of these works focus on specific families of low-dimensional dynamics: \cite{Blank} for rotations and twist maps, \cite{MR835874} and \cite{MR1709708} for the tent map of slope $\pm 2$, \cite{Boya-comp} for a piecewise expanding map of slope $\pm 3$, \cite{MR1258483} for $2x \mod 1$ (but for a random roundoff error model), \cite{MR1481914, MR1400185, MR1709708} for $1-|1-2x|^\ell$, \cite{MR1216205,MR1083721} for the Gauss map\dots{} These articles mainly focus on \emph{combinatorial} properties of discretizations: such discretizations are finite maps, hence their combinatorial properties are roughly determined by the family of lenghts of periodic orbits and the size of their basins of attraction. This focus on combinatorial properties seems to have been initiated in \cite{MR646380}.

In \cite{MR1678095}, after some illuminating general remarks, Lanford carries numerical simulations of the expanding map $x\mapsto 2x+x(1-x)/2$. Although mainly combinatorial, one of them computes the measure carried by the longest detected cycle. On these simulations, these measures seem close to the SRB measure (the discretization is taken in the sense of the double precision).

More recently, Galatolo, Nisoli and Rojas \cite[Sections 6, 7 and 8]{MR3187634} conducted numerical experiments on circle piecewise expanding maps (one example of which with a point with derivative 1) from an ergodic viewpoint. The difference with our study is that they consider Birkhoff averages of Dirac measures instead of the Lebesgue measure; their study is less extensive than ours but they still observe interesting behaviours of artefacts generated by roundoff errors. Their conclusion is that ``\emph{These experiments show that, in general, using floating point arithmetics to compute Birkhoff averages and invariant measures should not be considered reliable, not because of truncation and rounding errors, but rather because the dynamics of the discretised map does not mirror the generic dynamic of the real map.}''

There are few theoretical nontrivial results about the relations between discretizations and ergodic properties. In \cite{Gora-why}, G{\'o}ra and Boyarsky get some theoretical results of convergence of the discretizations' asymptotic measures $\mu_N$ (see \eqref{EqDefMuN}) towards $SRB$, under the hypothesis that there exists large orbits for the discretization (of size $\ge \alpha N$ for a fixed $\alpha>0$ and any $N$ large enough). They check that this hypothesis holds for some piecewise linear maps of slope that are power of $3$. However, \cite[Theorem 33]{MR3919917} (see Theorem~\ref{TauxExpand}) shows that this hypothesis is not satisfied for generic dynamics\dots

In his PhD thesis \cite{Flocker}, Flockermann carries numerical simulations of maps which are similar to those considered in the present paper. However, these are mainly combinatorial: as in \cite{MR1678095}, the only ergodic properties considered deal with the measure carried by some periodic orbits of the discretization. In this thesis the author also proves theoretical results about distribution of roundoff errors for generic circle expanding maps.

These results were obtained independently by Vladimirov in \cite{vladimirov2015quantized} (further works based on this grounding article were published in \cite{MR1661120, MR1769577, MR1956409, MR1894464}). In this article, the author founds a solid theoretical basis about the discretizations' behaviour, which revals more powerful than Flockermann's approach: in addition to the equidistribution of roundoff errors, Vladimirov gets Theorem~\ref{TauxExpand}, and some functional central limit theorem, which was published with Vivaldi in \cite{MR2031150}. Early apparitions of this kind of ideas can be found in the work of Voevodin \cite{MR219227}.

Part of these results were rediscovered independently (a second time!) by the first author in \cite{MR3919917}; this work also contains the case of diffeomorphisms and measure-preserving diffeomorphisms, it is based on an approach which, although a bit different, is quite similar to the one of \cite{vladimirov2015quantized}.

\section{Preliminaries}

\subsection{Distance on measures}

We will denote by $D^r(\Sp^1)$ the set of maps $f : \Sp^1\to\Sp^1$ that are $C^r$ and expanding (meaning that $f'(x)>1$ for any $x\in\Sp^1)$.

In the first paper of this series \cite{paper1}, we give an asymptotics of the distance between the measures $f_*^k(\Leb)$ and $(f_N)_*^k(\Leb)$, for $k$ fixed and $N$ going to infinity. This distance is measured by what we call the \emph{Cramér distance}: if $\mu$ and $\nu$ are two probability measures  on $\Sp^1$, and $F$ and $G$ are their respective repartition functions defined from the starting point 0, we set $H=F-G$ and
\begin{equation}\label{DefDisc}
\Disc(\mu,\nu) = \left(\min_{c\in\R} \int_0^1 (H-c)^2\right)^{1/2} = \left(\int_0^1 \left(H(x) - \Big(\int_0^1 H\Big)\right)^2 \ud x\right)^{1/2}.
\end{equation}
This is a distance spanning the weak-* topology on measures. For more details, see \cite{paper1}. The following lemma is straightforward.

\begin{lemma}\label{lemDiscrepDensit}
If $\mu$ and $\nu$ are two absolutely continuous probability measures on $\Sp^1$ with respective densities with respect to Lebesgue measure $f$ and $g$, then 
\[\Disc(\mu,\nu)\le \|f-g\|_1 \le \|f-g\|_1 \le \|f-g\|_\infty.\]
\end{lemma}

\begin{proof}
In this case we have, for any $x\in[0,1]$,
\[F(x) = \int_0^x f(t)\ud t\]
and the same for $G$, so
\[|H(x)| \le \int_0^1 |f-g|(t)\ud t = \|f-g\|_1, \]
and 
\[ \Disc(\mu,\nu) \le \left(\int_0^1 H^2\right)^{1/2} \le \|f-g\|_1 \le \|f-g\|_2 \le \|f-g\|_\infty.\]
\end{proof}

We will also use the \emph{Ruelle-Perron-Frobenius} (RPF) operator, defined on observables $\phi: \Sp^1\to\C$ by
\begin{equation}\label{EqDefRPF}
L_f \phi : y \mapsto \sum_{f(x)=y} \frac{\phi(x)}{f'(x)}.
\end{equation}
Note that if $\phi$ is the density of an absolutely continuous measure $\mu$, then $L_f\phi$ is the density of $f_*\mu$.

\subsection{The maps used in the simulations}\label{SecMaps}

In our numerical studies, we will consider the following maps:
\begin{align}
f_{c_1,c_2,k}: \Sp^1 & \longrightarrow \Sp^1 \label{EqDefDyna}\\
x & \longmapsto 2x+ c_1\sin(2\pi x) + c_2\sin(4\pi x) + k,\nonumber
\end{align}
with $c_1,c_2,k\in\R$ three parameters, with $c_1,c_2$ chosen such that the map $f$ is expanding (which is true if $2\pi |c_1|+4\pi |c_2| < 1$). 

In most of the simulations, we will take $c_1 = c_1^0 = 0.0531647$, $c_2 = c_2^0 = 0.03932758$ and $k=0.347$ (and in this case the minimum of $f_{c_1,c_2,k}'$ is bigger than $1.17$). See Figure~\ref{FigMap} for a graph of this map.

For some simulations, we will consider small perturbations of this system, by choosing $c_1 = c_1^0 + 0.001p_1$ and $c_2 = c_2^0 + 0.001p_2$, for
\[p_1,p_2\in \left\{ -1,-\frac12,0,\frac12,1\right\}.\]

\subsection{The code for the experiments}
 
The code we used for experiment is based on the \emph{Python} project
\emph{CompInvMeas-Python} \cite{compinvmeas} developed as an initiative to unify the
approach to the computation of invariant measures explained in
\cite{galatolo2014elementary} using \emph{SageMath} \cite{sagemath}, the framework and
related further developments will be fully described in the article in preparation
\cite{generalframework}.  The code used from the project was forked from an older version
and contains facilities to work with dynamical systems, to compute the Perron-Frobenius
operator of an expanding dynamical system and to retrieve the corresponding numeric fixed
point.

When the dynamics is expanding with a factor $\geq 2$ it is possible to certify the error
in the approximation of the SRB measure, estimating independently the numerical error
which occurred while computing the numeric fixed point, and the mathematical error
occurring representing the transfer operator with a finite-dimensional linear
operator. While relevant, this estimation is very pessimistic therefore we used the result
of the SRB measure in our experiments without adding the error coming from the rigorous
error estimation, as it would have hidden the error coming from the spatial
discretization.

Our experiments have been conducted in a notebook using the above facilities as a \emph{Python}
library, plus additionally a few support Python files offering facilities of our
experiments. Such facilities are to compute different spatial discretizations of
a dynamical system, to compute measure distances (Cramér and Wasserstein distance),
and convenience functions to save intermediate results to a database in order to be able
to interrupt experiments and resume them later. An end-to-end run of all the experiments
takes several days and will create temporary files of roughly 500Gb.

The notebook with all the support file and the instructions to repeat the experiments is
available on the url: \url{https://github.com/maurimo/DiscretizedDynSys}

\section{Short term behaviour}\label{SecShort}

\subsection{Theoretical result for Cramér distance}

In \cite{paper1}, we get an asymptotics for the map \eqref{EqStudiedQuant} for a generic expanding map.

\begin{theoreme}[\cite{paper1}]\label{MainTheo}
Let $r\ge 1$, $f$ a generic $C^r$ expanding map of the circle $\Sp^1$, and $k\in\N$. Then
\begin{equation}\label{EqIntMieuxIntro}
\lim_{N\to +\infty} N^2 \Disc\Big((f_N^k)_*(\Leb_N),\, f^k_*(\Leb)\Big)^2
= \frac{1}{12} + \frac{1}{12} \sum_{m=0}^{k-1} \big\langle D(f^{k-m}), (L_f^m 1)^2 \big\rangle,
\end{equation}
where $\langle\cdot,\cdot\rangle$ stands for the $L^2$ scalar product, $L_f$ is the RPF transfer operator defined by \eqref{EqDefRPF}, $D$ is the derivative and $f^{k-m}$ is the $(k-m)$-th iterate of $f$.
\end{theoreme}

\begin{figure}[ht]  
\centering \includegraphics[width=.6\linewidth]{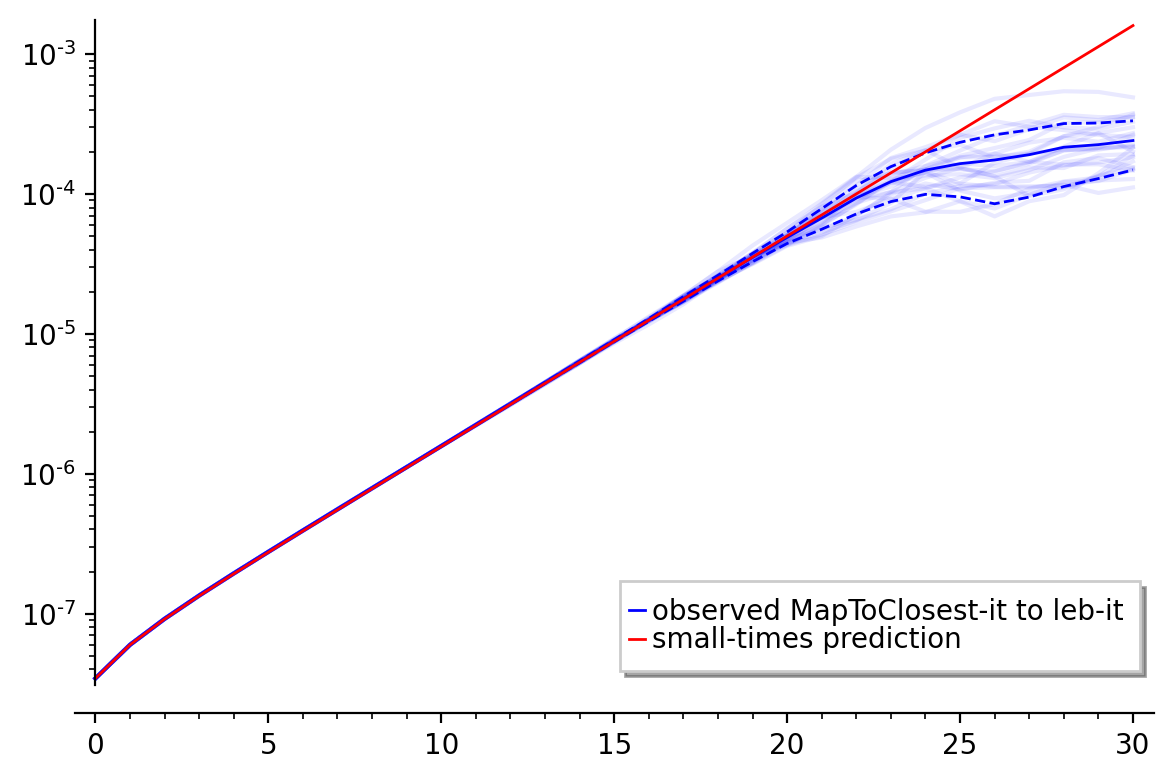}
\caption{\label{FigSmallTimes} Distance \eqref{EqStudiedQuant} depending on time $k$, for $N=2^{23}$, and for 25 perturbations of the map \eqref{EqDefDyna} described in Section~\ref{SecMaps} (in logarithmic vertical scale). The different curves for the different perturbations are in light blue, the blue curve represents the mean and the dashed curves the mean $\pm$ the standard deviation. The red curve is the theoretical prediction given by Theorem~\ref{MainTheo}, and computed with the help of the RPF operator \eqref{EqDefRPF} for which we have a fast approximation algorithm.}
\end{figure}

The aim of this section is to explore numerically the validity in practice of such results: the speed of convergence cannot be specified in the proof of Theorem~\ref{MainTheo} (it is hidden in the ``generic'' term).

\begin{figure}[ht]
\begin{minipage}{0.49\linewidth}   
\centering \includegraphics[width=\linewidth]{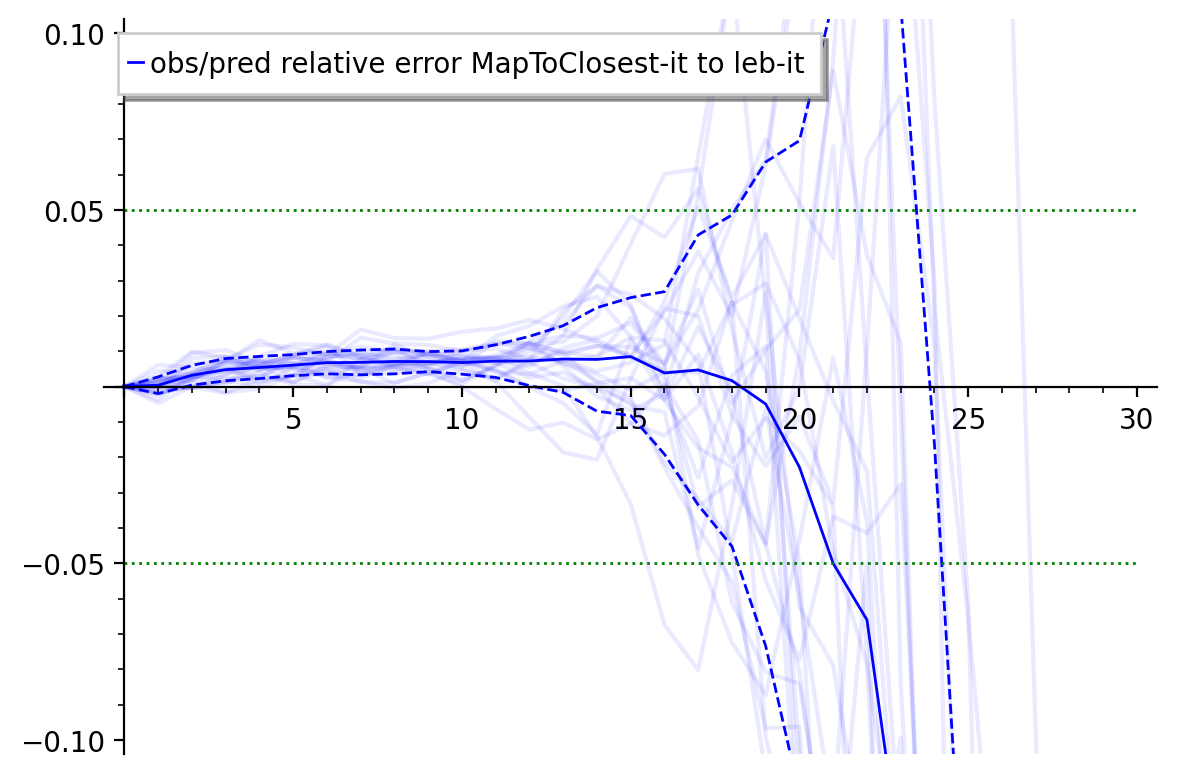}
\end{minipage}
\hfill
\begin{minipage}{.49\linewidth}
\centering \includegraphics[width=\linewidth]{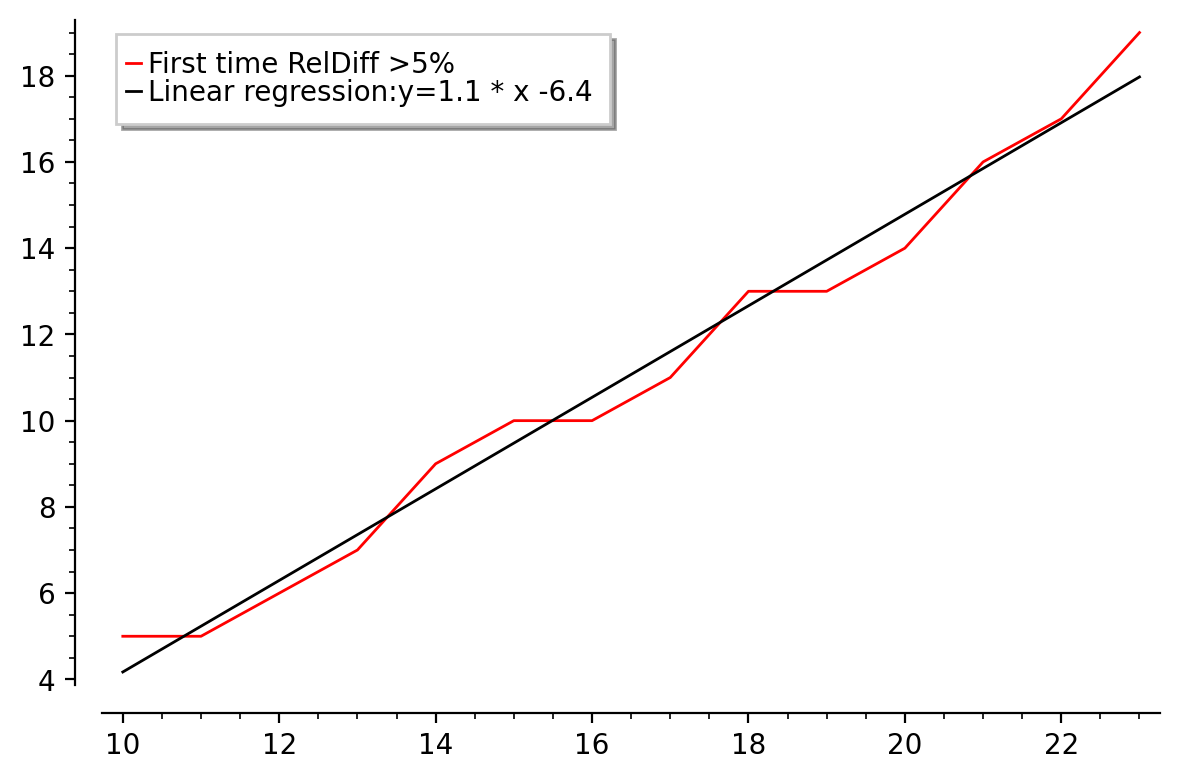}
\end{minipage}
\caption{\label{FigSmallTimes2} Left: relative difference between the theoretical prediction of Theorem~\ref{MainTheo} and the actual distances $\Disc$ on the examples for $N=2^{23}$. Right: the $x$-axis represents $k=\log_2(N)$, the $y$-axis represents the first time for which the mean of the relative difference $\pm$ standard deviation is bigger than $5\%$ (i.e. the first time one of the left blue dashed curves meets one of the green dotted lines) depending on $k=\log_2(N)$.}
\end{figure}

As can be seen on Figure~\ref{FigSmallTimes}, the theoretical prediction is quite good up to time 20 for $N=2^{23}$. More precisely, as can be seen on the left of Figure~\ref{FigSmallTimes2}, the theoretical prediction stops being relevant from time $\simeq 18$. Note that for $N=2^{23}$, one has $\log_2(N)=23$; this behaviour of the time until when the theoretical prediction is accurate typically logarithmic in $N$ is strongly suggested by Figure~\ref{FigSmallTimes2}, right. In can be explained heuristically in the following way: for $N=2^{23}$, as the derivative of the map $f$ is everywhere close to 2, 20 is more or less the time needed for the iterations of a grid domain $[i/N, (i+1)/N]$ to become macroscopically visible.

\begin{moral}
\textsf{In practice, Theorem~\ref{MainTheo} is valid until times logarithmic in $N$.}
\end{moral}

\subsection{Rate of injectivity}\label{SecRat}

Another quantity for which we have theoretical results is the \emph{rate of injectivity}. It is defined as
\[\tau^k(f_N) = \frac{\card\big((f_N)^k(E_N)\big)}{\card(E_N)}.\]
This quantity (and the one studied in the next subsection) will be used in the study of the medium term behaviour of discretizations. In this subsection and the next one, we will:
\begin{itemize}
\item state theoretical results for these quantities, which will be proved for generic maps and small number of iterations;
\item observe experimentally whether these theoretical results stay true in the short or medium term.
\end{itemize}

\begin{figure}[ht]
\begin{center}
\begin{tikzpicture}[scale=1.1]
\node (O) at (0,0){$y$};
\node (A) at (3,1){$x_{(1)}$};
\node (B) at (3,-1){$x_{(2)}$};
\node (C) at (6,1.5){$x_{(1,1)}$};
\node (D) at (6,.5){$x_{(1,2)}$};
\node (E) at (6,-.5){$x_{(2,1)}$};
\node (F) at (6,-1.5){$x_{(2,2)}$};
\draw (O) -- (A) node[sloped, midway, above]{$1/f'(x_{(1)})$};
\draw (A) -- (C) node[sloped, midway, above]{$1/f'(x_{(1,1)})$};
\draw (A) -- (D) node[sloped, midway, below]{$1/f'(x_{(1,2)})$};
\draw (O) -- (B) node[sloped, midway, below]{$1/f'(x_{(2)})$};
\draw (B) -- (E) node[sloped, midway, above]{$1/f'(x_{(2,1)})$};
\draw (B) -- (F) node[sloped, midway, below]{$1/f'(x_{(2,2)})$};
\end{tikzpicture}
\caption[Probability tree associated to the preimages of $y$]{The probability tree $T_k(y)$ associated to the preimages of $y$, for $k=2$ and $d=2$. We have $f(x_{(1,1)}) = f(x_{(1,2)}) = x_{(1)}$, etc.}\label{ProbTree}
\end{center}
\end{figure}
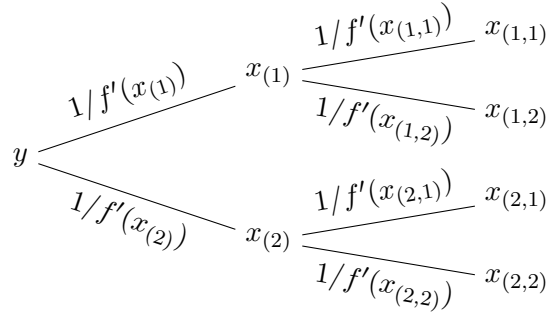

Before recalling the result of \cite{MR3919917}, let us introduce some notations. Given an expanding map $f$ of $\Sp^1$ of degree $d$, the set of time-$k$ preimages of a point $y\in\Sp^1$ has a structure of complete $d$-ary tree, whose vertices are the points $x\in f^{-m}(y)$ for $0\le m \le k$, and the edges are of the form $(x,f(x))$. One labels each edge $(x,f(x))$ of this tree by the number $1/f'(x)$, and denote by $T_k(y)$ the resulting labelled graph (see Figure~\ref{ProbTree}). 

We call \emph{random graph associated to $f$ at $y$} the random subgraph $G_k(y)$ of $T_k(y)$, such that the laws of appearance of the edges $(x,f(x))$ in $G_k(y)$ are independent Bernoulli laws of parameter $1/f'(x)$. In other words, $G_k(y)$ is obtained from $T_k(y)$ by erasing independently each vertex of $T_k(y)$ with probability $1-1/f'(x)$.

We define the \emph{mean density} $\overline D_k(y)$ as the probability that in $G_k(y)$, there is at least one path linking the root to a leaf.

The following is a restatement of \cite[Theorem 33]{MR3919917} (see also \cite{vladimirov2015quantized}).

\begin{theoreme}\label{TauxExpand}
Let $r\ge 1$, $f$ a generic element of $\mathcal D^r(\Sp^1)$ and $k\in\N$. Then,
\begin{equation}\label{EqIntMieux}
\lim_{N\to +\infty} \tau^k(f_N) = \int_{\Sp^1} \overline D_k(y) \ud \Leb(y).
\end{equation}
\end{theoreme}

As a byproduct of the proof of this theorem (and in particular Lemma~34 of \cite{MR3919917}), we get the following local convergence result (see also \cite{vladimirov2015quantized}).

\begin{prop}\label{PropLocDistrib0}
For any $r\ge 1$, for a generic expanding map $f\in\D^r(\Sp^1)$ and for almost every point $y$, one has
\[\overline D_k(y) = \lim_{\substack{N,R\to +\infty \\ R/N\to 0}} \frac{1}{2R} \card\Big\{\frac{i}{N}\in [y-R/N, y+R/N]\ \big|\ \frac{i}{N}\in f_N^{k}(E_N)\Big\}.\]
\end{prop}

Note that a first step towards the proof of this theorem was realized in the unpublished thesis \cite{Flocker}.

The idea behind this theorem is the following. Assume for simplicity that $d=2$, take some point $y\in \Sp^1$, and denote its preimages by $f$ by $x_0$ and $x_1$. Then, in the neighbourhood of $y$, the set $f_N(E_N)$ looks like the discretization\footnote{Here, ``discretization'' stands for the projection on the nearest element of $\Z$, i.e. the image under the projection $\R\to\Z$ on the nearest integer.} of the set $f'(x_0)\Z \cup f'(x_1)\Z$. But, still in the neighbourhood of $y$, the ``probability'' for a point $z\in E_N$ to be in the discretization of the set $f'(x_0)\Z$ is equal to $1/f'(x_0)$. One of the steps in the proof of Theorem~\ref{TauxExpand} is to show that the probabilities coming from the different branches are independent: the probabilities for a point $z\in E_N$ to be in the discretizations of the sets $f'(x_i)\Z$ are independent.
\bigskip

The following lemma gives a practical way to compute the percolation probability $\overline D_k(y)$, in terms of the transfer operator $L_f$ associated to $f$ (for which we have a fast and reliable algorithm).  

\begin{lemma}\label{LemDk}
For any $y\in \Sp^1$,
\[\overline D_{k+1}(y) = 1-\prod_{x\in f^{-1}(y)}\left(1-\frac{\overline D_{k}(x)}{f'(x)}\right).\]
In particular, if the degree of $f$ satisfies $d=2$, denoting $f^{-1}(y) = \{x_0,x_1\}$, one has
\begin{align*}
\overline D_{k+1}(y) & = \frac{\overline D_{k}(x_0)}{f'(x_0)} + \frac{\overline D_{k}(x_1)}{f'(x_1)} - \frac{\overline D_{k}(x_0)\overline D_{k}(x_1)}{f'(x_0)f'(x_1)}\\
& = L_f(\overline D_{k}) - \frac12\left( \big(L_f(\overline D_{k})\big)^2 - L_f\Big( \overline D_{k}^2/f'\Big)\right).
\end{align*}
\end{lemma}

It is possible to get similar formulae for bigger $d$ by using Vieta's formulas.

\begin{proof}
The first formula comes directly from the definition. The second one is a simple computation.
\end{proof}
\medskip

\begin{figure}[ht]
\includegraphics[width=.49\linewidth]{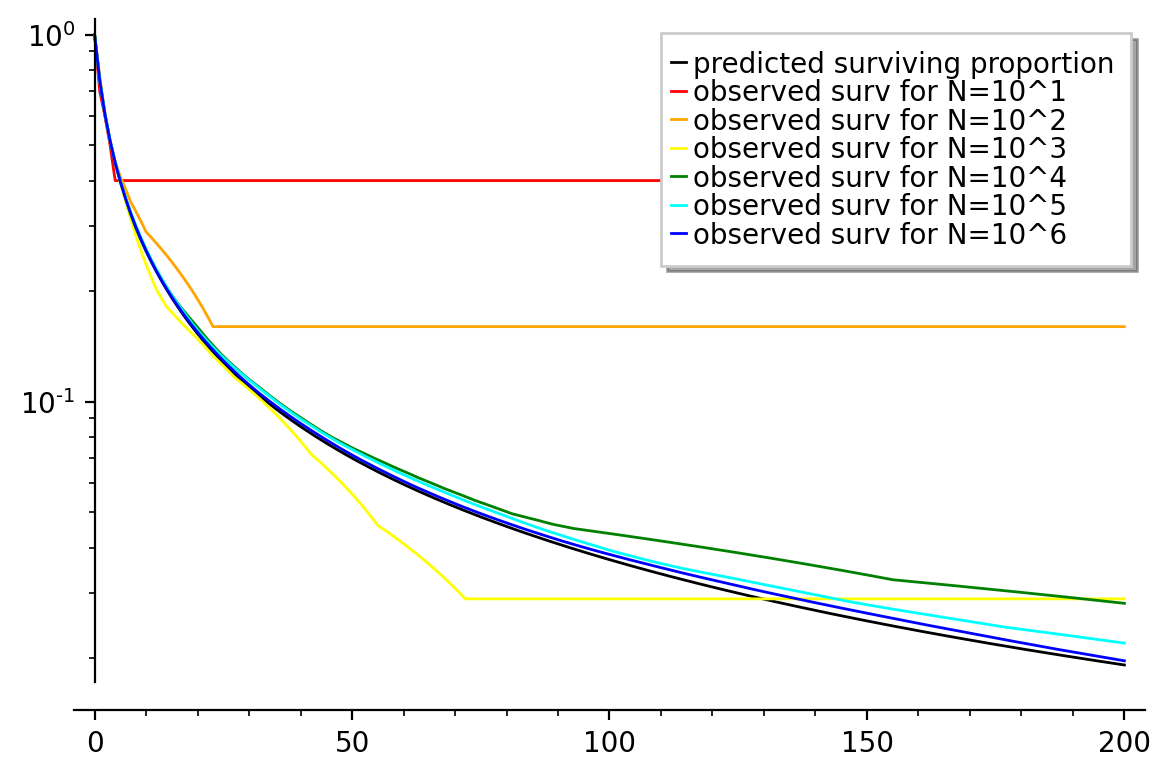}
\hfill
\includegraphics[width=.49\linewidth]{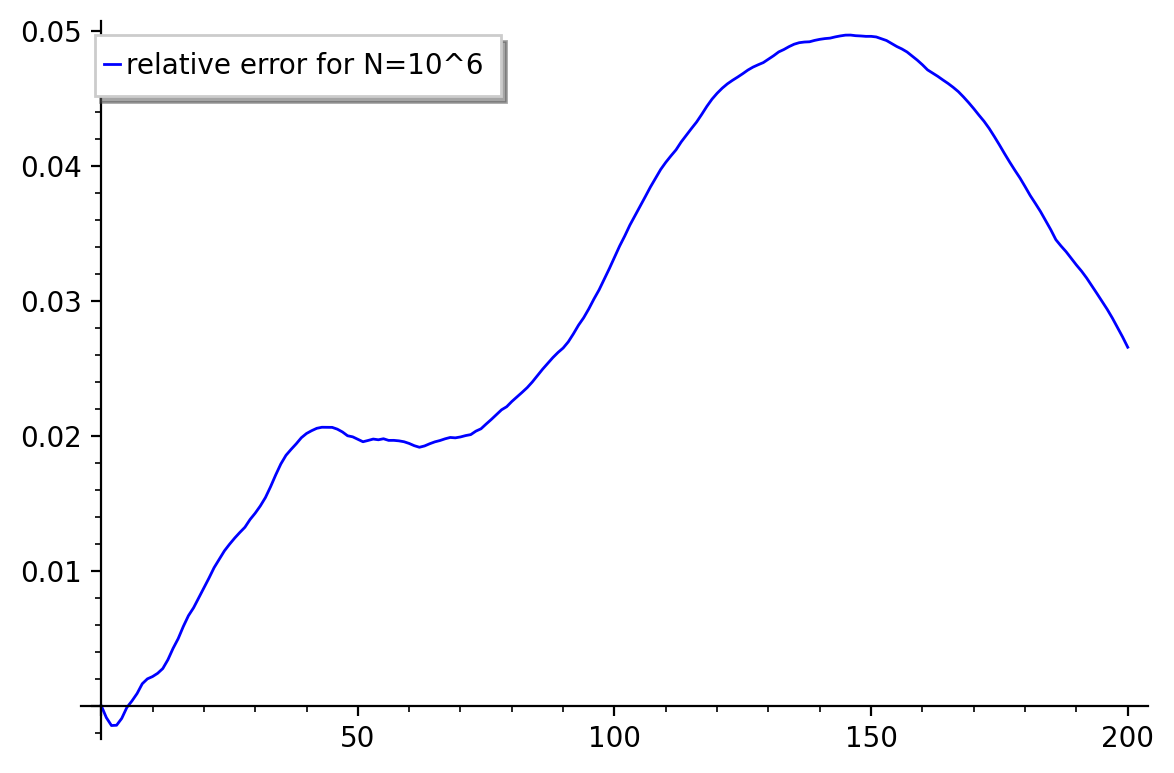}
\caption{\label{RateInjecTps} Left: rate of injectivity of $\tau^k(f_N)$ (the coloured curves) depending on $k$, for 6 different values of $N$: $10^1, 10^2, 10^3, 10^4, 10^5$ and $10^6$. We also represent the theoretical value (given by Theorem \ref{TauxExpand}, black curve) depending on $k$ in logaritmic scale.
Right: relative difference between these two quantities depending on $k$, for the biggest $N$ ($N = 10^6$).}
\end{figure}

Figure \ref{RateInjecTps} compares the theoretical rate of injectivity (the one given by Theorem \ref{TauxExpand} and computed with the help of Lemma~\ref{LemDk}) with the actual rate of injectivity of a discretization $f_N$, depending on the time $k$.

\begin{moral}
\textsf{The predictions of Theorem \ref{TauxExpand} are really good during a quite long time (less than 5\% up to time 200).}
\end{moral}

This contrasts with what happens for predictions of Theorem \ref{MainTheo} (see Figure~\ref{FigSmallTimes}), that become inaccurate in times typically logarithmic in $N$. We do not have any clue why these predictions stay accurate for such a long time.

\subsection{Local distribution of preimages}\label{SeclocDistrib}

We pursue the study of the rate of injectivity by focusing on more precise quantities. We will look at the \emph{distribution} of the number of preimages of a point of the grid. To do that, we set $a_m$ (which depends on the point $y$ and the time $k$) the probability that in $G_k(y)$, there are exactly $m$ paths linking the root with the leaves. Denote by
\begin{equation}\label{EqDefAi}
\overline P_k(y) = \sum_{m\ge 0} a_m X^m
\end{equation}
the associated generating series.

Of course, the polynomial $\overline P_k(y)$ is of degree at most $d^k$ and satisfies $\overline P_k(y)(1)=1$ and $\overline P_k(y)(0) = 1-\overline D_{k}(y)$.

The proof of Theorem \ref{TauxExpand} links the limiting behaviour, for $N\to +\infty$, of the number of preimages of points of the grid $E_N$, with the random behaviour of the tree $G_k(y)$. Hence, as a byproduct of this proof, one gets the following result (see also \cite{vladimirov2015quantized}).

\begin{prop}\label{PropLocDistrib1}
For any $r\ge 1$, for a generic expanding map $f\in\D^r(\Sp^1)$ and for almost every point $y$, one has
\[a_m = \lim_{\substack{N,R\to +\infty \\ R/N\to 0}} \frac{1}{2R} \card\Big\{i/N\in [y-R/N, y+R/N]\ \big|\ \card\big(f_N^{-k}(i/N)\big)=m\Big\}.\]
\end{prop}

\begin{figure}[ht]
\begin{minipage}{0.48\linewidth}   
\centering \includegraphics[width=\linewidth]{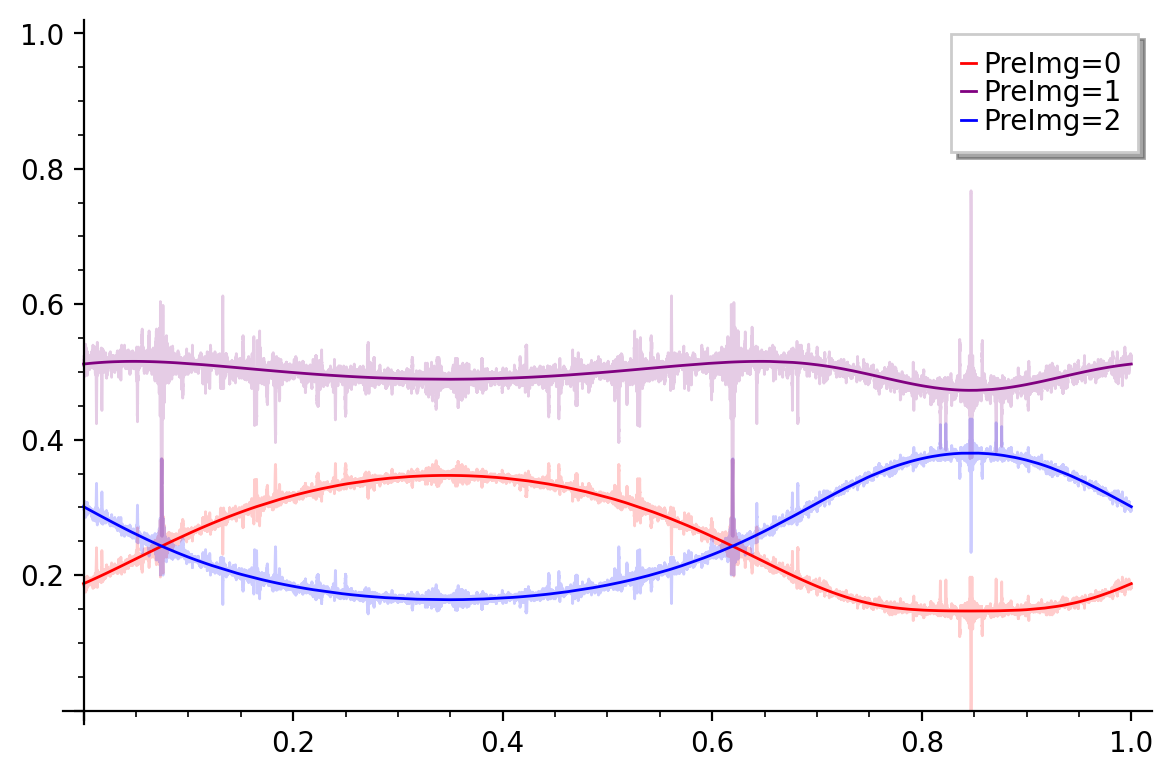}
\end{minipage}
\hfill
\begin{minipage}{.48\linewidth}
\centering \includegraphics[width=\linewidth]{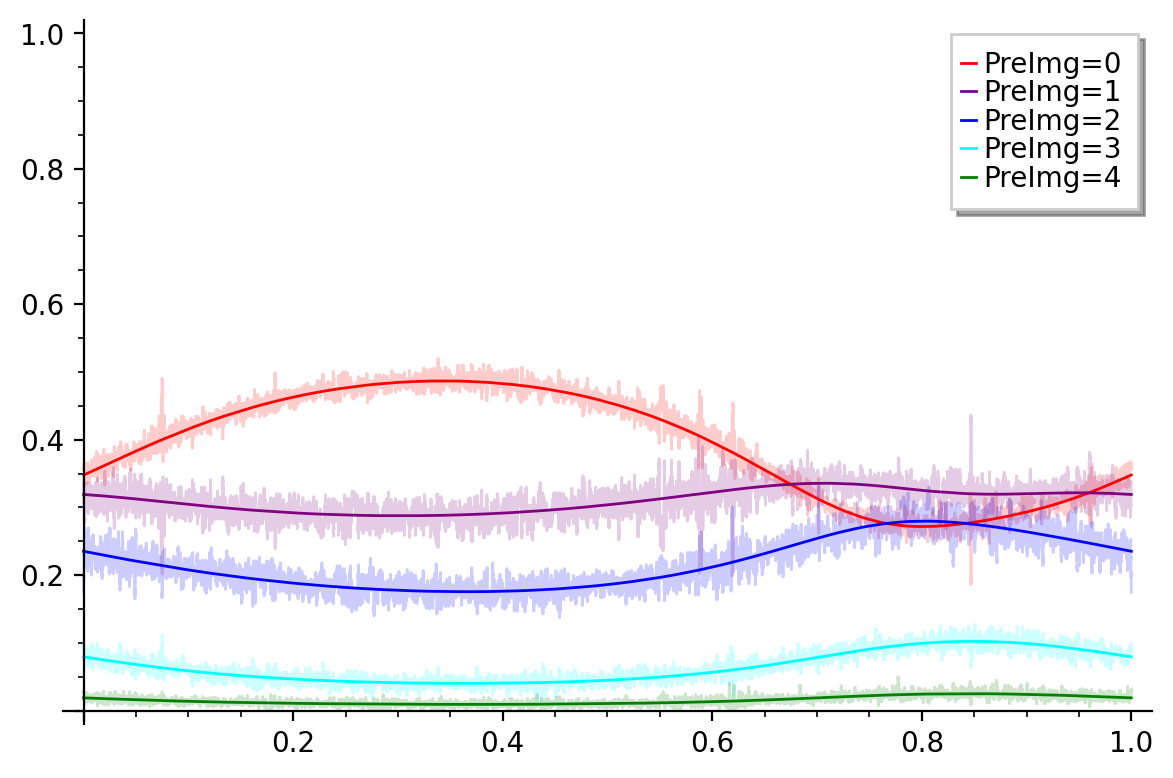}
\end{minipage}\\
\begin{minipage}{0.48\linewidth}   
\centering \includegraphics[width=\linewidth]{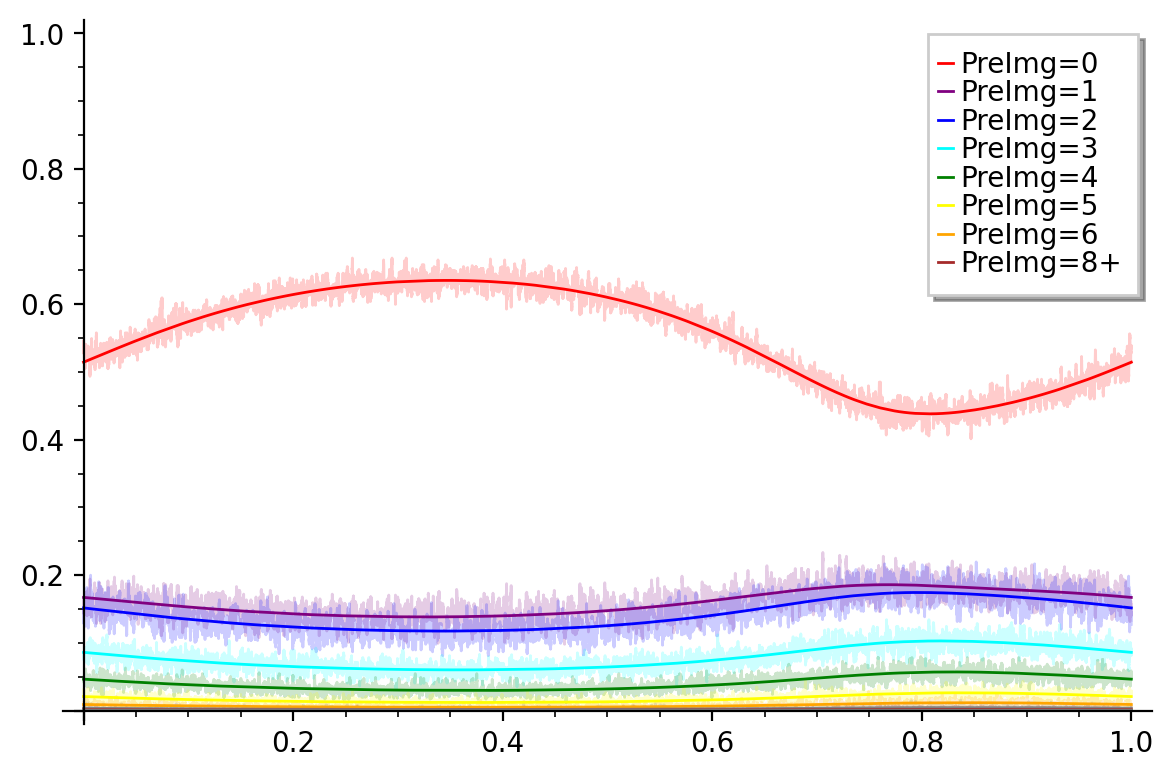}
\end{minipage}
\hfill
\begin{minipage}{.48\linewidth}
\centering \includegraphics[width=\linewidth]{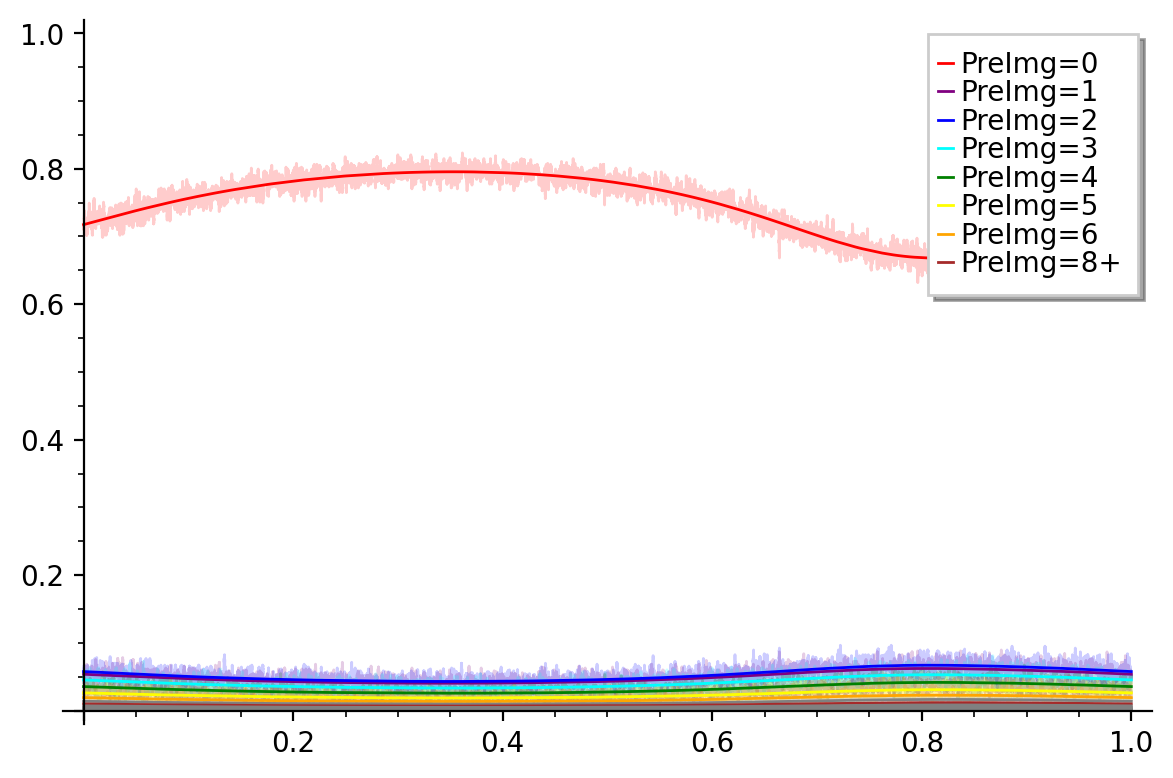}
\end{minipage}

\caption{\label{RateInjecDens}These curves represent the local densities of the number of preimages. The expected theoretical values given by Proposition~\ref{PropLocDistrib1} --- and computed with the help of Proposition~\ref{LemLocDistrib} --- are represented in full colors. The actual values for the map are in light colours, they represent the quantity \eqref{EqRateInjecDens} for $R=40$ and
$N = 10^6$. The different graphs correspond to different times: from left to right and top to bottom, $k=1, 2, 4$ and $10$.}
\end{figure}

In this case, the generating series formalism gives the following nice formula, that allows to compute the distributions $a_m$ by an iterating process involving the RPF operator $L_f$.

\begin{prop}\label{LemLocDistrib}
If $d=2$, then, denoting $f^{-1}(y) = \{x_0,x_1\}$,
\begin{align*}
\overline P_{k+1}(y) = & \frac{1}{f'(x_0)}\left(1-\frac{1}{f'(x_1)}\right) \overline P_{k}(x_0) + \frac{1}{f'(x_1)} \left(1-\frac{1}{f'(x_0)}\right)\overline P_{k}(x_1)\\
& + \frac{1}{f'(x_0)f'(x_1)} \overline P_{k}(x_0)\overline P_{k}(x_1)\\
& + \left(1-\frac{1}{f'(x_0)}\right)\left(1-\frac{1}{f'(x_1)}\right)
\end{align*}
Hence, denoting $\overline Q_k = \overline P_k-1$,
\[\overline Q_{k+1} = \frac12 \big(L_f(\overline Q_k)\big)^2+ L_f\left(\overline Q_k-\frac12 \frac{\overline Q_k^2}{f'}\right).\]
\end{prop}

\begin{figure}[ht]
\centering
\includegraphics[width=.55\linewidth]{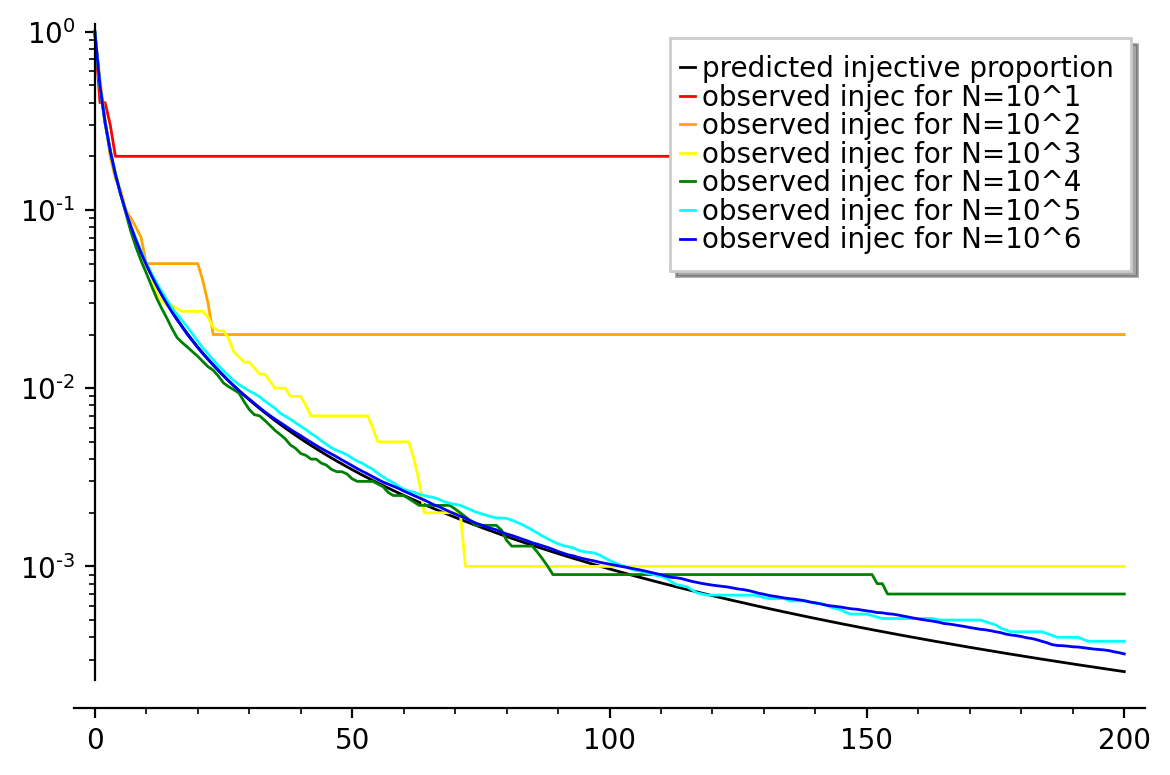}

\caption{\label{RateInjecTps2} Proportion of points of $E_N$ with one preimage under $f_N^k$ (the coloured curves) depending on $k$, for 6 different values of $N$: $10^1, 10^2, 10^3, 10^4, 10^5$ and $10^6$. We also represent the theoretical value (given by Proposition~\ref{PropLocDistrib1}, black curve) depending on $k$ in logaritmic scale.}
\end{figure}

\begin{proof}
Direct computation of the probabilities.
\end{proof}

This proposition gives a fast algorithm to compute the distributions $a_m$, as it allows to get it from iterations of $L_f$ for which we have a fast algorithm.

Figure \ref{RateInjecDens} gathers both theoretical and real values of the local densities of preimages. For the actual values, we represent the quantities, for different integers $m$, for a fixed discretization size $N$ and a fixed time $k$,
\begin{equation}\label{EqRateInjecDens}
\frac{1}{2R}\card\left\{ y\in E_N\cap [x-R/N, x+R/N]\ \middle\vert \ \card  f_N^{-k}(y) = m\right\}.
\end{equation}
As can be seen on this figure, the predictions of Proposition \ref{PropLocDistrib1} are quite accurate: the theoretical and actual curves match very well, up to time $10$, for $N=10^6$. In fact, these predictions stay accurate during quite a long time, as can be observed on Figure~\ref{RateInjecTps2}: for $N=10^5$ or $N=10^6$, there is no big difference between the observed values and the prediction up to time 200. In fact, it seems that the theoretical predictions of Figures~\ref{RateInjecTps} and \ref{RateInjecTps2} stay accurate until a time comparable to $\sqrt{N}$.

\begin{moral}
\textsf{The predictions of Theorem \ref{TauxExpand} and Propositions~\ref{PropLocDistrib0} and \ref{PropLocDistrib1} stay accurate until a time proportional to $\sqrt{N}$.}
\end{moral}

This suggests that the discretizations begin to deviate from these theoretical predictions when there is a significant proportion of points of $E_N$ that have fallen in a periodic orbit of $f_N$.

\section{Medium term behaviour}\label{SecMiddle}

In \cite{MR1678095}, Oscar E. Lanford proposes to study the dynamics of the maps $f_N^t$ in the regime $\log N \ll t \ll N$. Note that in the view of the discussion of Section \ref{SecAsym}, one might be tempted to replace the last condition by $\log t \ll \log N$, as sugested by Lanford himself (see also Figure~\ref{TimeToCycle} and the associated discussion). For the first condition $\log N \ll t$, it is also justified by the previous study of the short term behaviour, as well as \cite[Theorem 12.17]{Guih-These}.

For now, theoretical breakthroughs in Landford's regime seem out of reach, which motivates a numerical study of discretizations in this case.

Recall the different phenomena isolated in the introduction, that can explain why the action of discretizations $f_N$ on measures differ from the one of the initial map $f$.
\begin{enumerate}[label=($P_\arabic*$), ref=($P_\arabic*$)]
\item\label{P1} The iterates of points always belong to $E_N$.
\item\label{P2} Two points of $E_N$ having the same image by $f_N$ will have identical positive orbits.
\item\label{P3} The local shape around $y\in\Sp^1$ of the image $f_N(E_N)$ is very similar to the one of a linearization of $f$ around the points $f^{-1}(y)$, which is a model set (see \cite{paper1}).
\item\label{P4} Any point eventually falls in a periodic cycle.
\end{enumerate}

\subsection{Different discretization schemes}\label{ParagDifferent}

To understand what influences the evolution of the distance $\Disc$ between iterates of Lebesgue measure and iterates of the uniform measure under discretizations, we look at what happens when we change the definition of the discretized map. We will need the following notation: for $x\in E_N$, the integer $i_x$ is chosen such that $i_x/N \le f(x) < (i_x+1)/N$; it allows to set $\epsilon_x \in[0,1]$ such that $f(x) = (i_x+\epsilon_x)/N$.
\begin{itemize}
\item \texttt{MapToClosest}: This is the already defined discretization $f_N$ of the map, where $f_N(x)$ is defined as the point of $E_N$ closest to $f(x)$.
\item \texttt{OnceDecidedRandom}: $D_N^o(f) : E_N\to E_N$ is a random map, such that for each $x\in E_N$, the point $D_N^o(f)(x)$ is chosen \emph{once for all} and randomly (and independently) to be $i_x/N$ with probability $1-\epsilon_x$, and to be $(i_x+1)/N$ with probability $\epsilon_x$. Note that the iterations of two points $x,y\in E_N$ under $[D_N^o(f)]^2$ are independent iff $D_N^o(f)(x) \neq D_N^o(f)(y)$.
\item \texttt{StepwiseRandom}: $D_N^s(f) : E_N\to E_N$ is a random map quite similar to $D_N^o(f)$ (\texttt{OnceDecidedRandom}), such that for each $x\in E_N$ \emph{and at each iteration}, the point $D_N^s(f)(x)$ is chosen randomly (and independently) to be $i_x/N$ with probability $1-\epsilon_x$, and to be $(i_x+1)/N$ with probability $\epsilon_x$.
\item \texttt{PointsRandomOnGrid}: $D_N^g(f)$ acts independently on $N$-tuples of elements of $E_N$ as $D_N^s(f)$ (\texttt{StepwiseRandom}): 
\[D_N^g(f)(x_1,\dots,x_N) = \big(D_N^s(f)(x_1),\dots,D_N^s(f)(x_N) \big).\]
Of course, this gives the measure
\[\frac{1}{N}\sum_{i=0}^{N-1} \delta_{[D_N^s(f)]^k(i/N)}.\]
\item \texttt{PointsPerturbed}: $D_N^p(f)$ acts on $N$-tuples of elements of $\Sp^1$ as a random perturbation of $f$: let $\tilde f_N$ be the random map obtained from $f$ by post-composing with a uniform noise on the segment $[-1/(2N),1/(2N)]$ (\emph{i.e.} $\tilde f_N(x)$ is chosen randomly and uniformly in $[f(x)-1/(2N),f(x)+1/(2N)]$). Then 
\[D_N^p(f)(x_1,\dots,x_N) = \big(\tilde f_N(x_1),\dots,\tilde f_N(x_N) \big).\]
\item \texttt{MapToCombination}: $D_N^c(f)$ acts only on the measures on $E_N$. It is affine, in the sense that for any convex combination $\mu = \sum_i\lambda_i\delta_{x_i}$, one has $D_N^c(f)(\mu) = \sum_i\lambda_i D_N^c(f)(\delta_{x_i})$. And $D_N^c(f)(\delta_x)$ is defined by
\[D_N^c(f)(\delta_x) = (1-\epsilon_x) \delta_{i_x/N} + \epsilon_x \delta_{(i_x+1)/N}.\]
\end{itemize}

Let us discuss the fundamental differences between these maps.

\texttt{MapToClosest} and \texttt{OnceDecidedRandom} have a quite similar definition, except that the first one is deterministic and the second one random. More precisely, there is the following difference between these maps: as $f$ is almost linear at a small scale, the image of $(i+1)/N$ under $f_N$ will depend \emph{deterministically} on $f(i/N)\mod N$, while the image of $(i+1)/N$ under $D_N^o(f)$ will depend only \emph{probabilistically} on $f(i/N)\mod N$. In other words, for $N$ large enough, $f_N(E_N)$ is locally almost (i.e. up to a set of arbitrarily small local density) a \emph{model set} (for a definition and a study of this property, see \cite{paper1} or \cite{MR3919917}), a property that $D_N^o(f)(E_N)$ does not possess. This will allow us to determine if the phenomenon \ref{P3} has a detectable effect on the evolution of the distance $\Disc$ \eqref{EqStudiedQuant} between $(f_N^k)_*(\Leb_N)$ and $f^k_*(\Leb)$.

The difference between \texttt{OnceDecidedRandom} and \texttt{StepwiseRandom} is that the second one is not autonomous. Hence, almost surely, orbits for this map will not be pre-periodic, while all orbits of \texttt{MapToClosest} and \texttt{OnceDecidedRandom} are. This will allow us to determine if the phenomenon \ref{P4} has an effect on the evolution of the distance \eqref{EqStudiedQuant}.

An important feature of the three previous maps (\texttt{MapToClosest}, \texttt{OnceDecidedRandom} and \texttt{StepwiseRandom}) is that orbits that merge then stay together forever. The map \texttt{PointsRandomOnGrid}, which besides is quite similar to \texttt{StepwiseRandom}, does not have this property. A priori, $D_N^g(f)(x) = D_N^g(f)(y)$ does not imply that $[D_N^g(f)]^2(x) = [D_N^g(f)]^2(y)$ a.s. This will allow us to determine if \ref{P2} affects the evolution of the distance \eqref{EqStudiedQuant}.

All the four previous maps are based on the discretization grid $E_N$. This is not the case for \texttt{PointsPerturbed}, which acts on $N$-tuples of points of the circle. It can be seen as a continuous counterpart of \texttt{PointsRandomOnGrid}. This will allow us to determine if the phenomenon \ref{P1} affects the evolution of \eqref{EqStudiedQuant}.

Finally, \texttt{MapToCombination} is the only discretization type which splits measures \footnote{By Perron-Frobenius theorem, the measures obtained from \texttt{MapToCombination} tend (when the time goes to infinity) to some measure depending on $N$. We do not know if these measures tend to $SRB$ when $N$ goes to infinity. It may be possible to prove it using the ideas of \cite{galatolo2014elementary}, by checking that the distance between both Perron-Frobenius and discretized (associated to \texttt{MapToCombination}) transfer operators are close relative to $BV_1$ distance. As this is not in the scope of this article, we do not investigate this question.}.

\begin{figure}[ht]
\includegraphics[width=.49\linewidth]{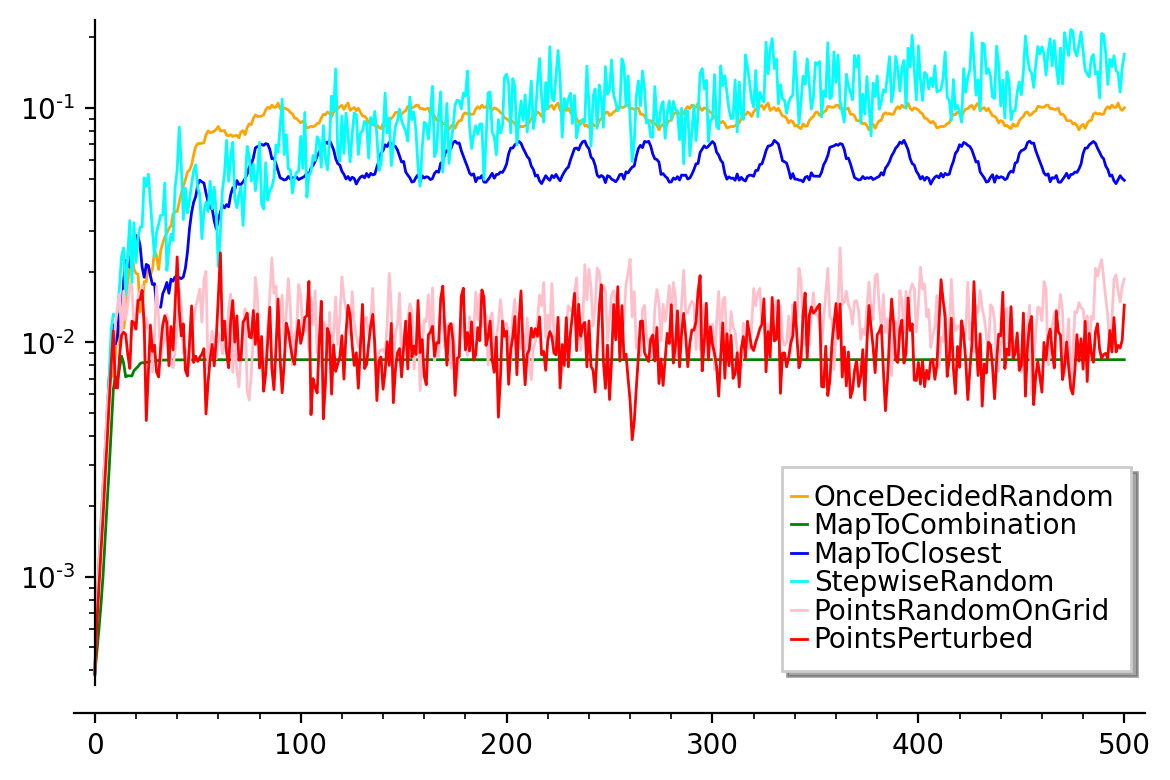}\hfill
\includegraphics[width=.49\linewidth]{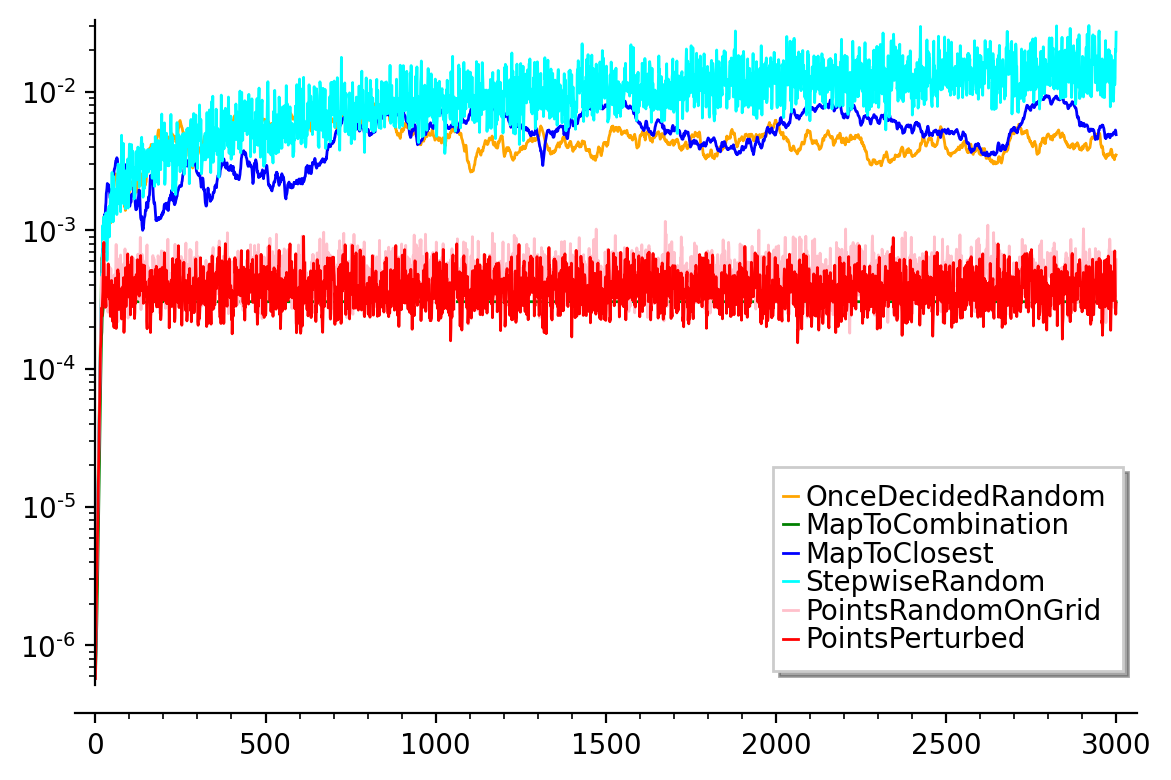}
\caption{\label{DiscrepancyIteratesDisc} Distance $\Disc$ between the measures $f^k_*(\Leb)$ and the images of $\Leb_N$ under $f_N^k$ and $(D_N^i(f))^k$ for $i\in\{o,s,g,p,c\}$ depending on the time $k$. On the left graphic, $N=750$ and on the right one, $N=50\,000$.}
\end{figure}

Figure~\ref{DiscrepancyIteratesDisc} shows the evolution of the distance $\Disc$ between the measures $f^k_*(\Leb)$ and the images of $\Leb_N$ under the iterates of the different discretization types of $f$: $f_N$, $D_N^o(f)$, $D_N^s(f)$, $D_N^g(f)$, $D_N^p(f)$ and $D_N^c(f)$.

On the left of this figure, where $N=750$, there is no intermediate regime for the discretizations $f_N$ and $D_N^o(f)$: from time $k\simeq 100$, the distance evolution becomes periodic. This can be explained by the fact that from this time, most of the grid's points have fallen in a few cycles of the discretization: one directly jumps from the small term behaviour, which is described quite well by Theorem~\ref{MainTheo}, to the periodic asymptotic regime. Indeed, in this case the limit time of short-term behaviour should be around $\log_2 N\simeq 10$, while the theoretical average time for orbits to cycle is $\sqrt{\pi N/2}\simeq 34$ (see Section~\ref{SecAsym}).

For $N=50\,000$ (Figure~\ref{DiscrepancyIteratesDisc}, right), this periodic asymptotic behaviour does not appear clearly for time $k\le 3\,000$: there is an actual medium term behaviour of discretizations.

Before studying further this intermediate regime, we first examine random point processes in the point of view of the distance $\Disc$ with the initial measure.

\subsection{Cramér distance between a measure and the random point process associated to it}\label{ParagPointProcess}

We have seen that on the simulations we made, the theoretical estimates on the local distributions of preimages (Paragraphs~\ref{SecRat} and~\ref{SeclocDistrib}) stay relevant in the middle term. Hence, they can be used to set a conjectural behaviour of the distance $\Disc$ in the middle term. Let us first introduce some definitions.

Let $\tilde \mu_1,\dots,\tilde \mu_M$ be probability measures on $\Sp^1$, with the density of $\tilde\mu_m$ being given by the map $a_m/\int_{Sp^1}a_m$ (recall that $a_m$ was defined in \eqref{EqDefAi} and described in Proposition~\ref{PropLocDistrib1}).
Let $K\in\N$, $(m_i)_{1\le i\le K}\in\{1,\dots,M\}^K$, and set $N = \sum_{i=1}^K m_i$.
Let also $\nu_\mu$ be the random probability measure defined by
\[\nu_\mu = \frac{1}{N} \sum_{i=1}^{K}m_i\delta_{p_i},\]
where each point $p_i$ is chosen independently in $\Sp^1$, according to the measure $\tilde \mu_{m_i}$.
Suppose that for any $m$,
\[\card\{i\mid m_i=m\}\simeq \int_{\Sp^1} a_m\]
(meaning that we have this asymptotics when $N$ goes to infinity).

\begin{conj}
In the regime $\log N \ll t$ and $\log t \ll \log N$, the distance\eqref{EqStudiedQuant} is close  to the expected value of the Cramér distance between the SRB measure and $\nu_\mu$.
\end{conj}

As a first step, before getting to the numerical study, we compute the expected value of the square of distance $\Disc$ between the SRB measure and $\nu_\mu$.

Let $\mu$ be a probability measure on $\Sp^1$. We identify $\Sp^1$ with $[0,1]$, and define $f$ as\footnote{There is a conflict of notations with the dynamics $f$, we hope that which one is used is clear from the context.} the cumulative distribution function of $\mu$ minus its average (so that $\int_0^1 f = 0$). Let also $F$ be the primitive function of $f$ such that:
\[F(x) = \int_0^x f(t) \ud t.\]
Remark that $F(0) = F(1) = 0$, and that $F\le 0$. Finally, given $p\in\Sp^1$, one defines $g_p(x) = \chi_{[p,1]}(x)-(1-p)$ to be the cumulative-minus-average distribution function of $\delta_p$

The following theorem gives the expectation of the (square of the) distance $\Disc$ between a measure $\mu$ and a point process associated to it.

\begin{theoreme}\label{TheoExpectDisc}
Let $\mu$ and $\tilde \mu_1,\dots,\tilde \mu_M$ be probability measures on $\Sp^1$ with respective repartition functions $f$ and $\tilde f_1,\dots,\tilde f_M$.

Let $(m_i)_{1\le i\le K}\in\{1,\dots,M\}^K$, and set $N = \sum_{i=1}^K m_i$.
Let also $\nu_\mu$ be the random probability measure defined by
\[\nu_\mu = \frac{1}{N} \sum_{i=1}^{K}m_i\delta_{p_i},\]
where the points $p_i$ are chosen independently in $\Sp^1$, each one with distribution $\tilde \mu_{m_i}$.
Then
\begin{align}\label{EqExpectDisc2}
\mathbb{E}\left[\Disc(\mu,\nu_\mu)^2\right]
& =  \int_0^1\left(f-\sum_{i=1}^K \frac{m_i}{N}\tilde f_{m_i} \right)^2 - \sum_{i=1}^K \frac{m_i^2}{N^2} \int_0^1\left(\tilde f_{m_i}^2 + 2\tilde F_{m_i}\right)\nonumber\\
& =  \Disc\left(f,\,\sum_{i=1}^K \frac{m_i}{N}\tilde f_{m_i} \right)^2 + \sum_{i=1}^K \frac{m_i^2}{N^2} \left(\frac{1}{12} - \Disc\left(\tilde f_{m_i},\, \Leb\right)^2\right).
\end{align}
\end{theoreme}

The proof of this theorem can be found in appendix.

The setting of this theorem will be applied in the case where $\mu$ can be written as $\mu = \sum_{m=0}^M m_i \mu_i$ (hence the $\mu_i$ are not probability measures), and the measures $\tilde \mu_i$ are the normalizations of the $\mu_i$, with the $m_i$ being chosen such that the normalization factors are close to $m_i/N$ (see Remark~\ref{RemExpectDisc}).

Some similar estimations were obtained for the Wasserstein distance $W_1$ in \cite{Bobkov2019OnedimensionalEM}, but the authors only manage to get bounds ant not exact values for the expected value.

\begin{rem}\label{RemExpectDisc}
Note that $N$ being fixed, one can choose the family $m_i$ such that the first term is of order $C/N^2$, while the second one is typically of order $1/N$.

Indeed, suppose that the cumulative distribution function $f$ of $\mu$ satisfies $f = \sum_{m=1}^K \lambda_m \tilde f_m$, and for $p\in\N$, let $(m_i)$ such that 
\[\card\{i\mid m_i = m\} = \lfloor p\lambda_m\rfloor.\]
Note that 
\[\frac{p\lambda_m-1}{\sum_{n=1}^K( p\lambda_n+1)} \le \frac{\lfloor p\lambda_m\rfloor}{\sum_{n=1}^K\lfloor p\lambda_n\rfloor} \le \frac{p\lambda_m+1}{\sum_{n=1}^K( p\lambda_n-1)},\]
so 
\[- \frac{1+K}{p+K} \le \frac{\lfloor p\lambda_m\rfloor}{\sum_{n=1}^K\lfloor p\lambda_n\rfloor} - \lambda_m \le \frac{1+K}{p-K},\]
and
\[\left\|f-\sum_{m=1}^K\frac{\lfloor p\lambda_m\rfloor}{\sum_{n=1}^K\lfloor p\lambda_n\rfloor} \tilde f_m \right\|_2 \le \sum_{m=1}^K \left|\frac{\lfloor p\lambda_m\rfloor}{\sum_{n=1}^K\lfloor p\lambda_n\rfloor} - \lambda_m\right|\|\tilde f_m\|_2 \le K\frac{1+K}{p-K}.\]
Hence,
\[\left\|f-\sum_{m=1}^K\frac{\lfloor p\lambda_m\rfloor}{\sum_{n=1}^K\lfloor p\lambda_n\rfloor} \tilde f_m \right\|_2^{2} = O\left(\frac{1}{p^2}\right),\]
which gives a distance of the order of $1/p^2$ by Lemma~\ref{lemDiscrepDensit}.

On the other hand, 
\[\sum_{m_i=m} \frac{m_i^2}{N^2}
= m^2 \frac{\lfloor p\lambda_m\rfloor}{\sum_{n=1}^K\lfloor p\lambda_n\rfloor} \frac{1}{\sum_{n=1}^K\lfloor p\lambda_n\rfloor}
\ge \frac{m^2\lambda_m}{p+K}\]
for any $p$ large enough.
Hence, in this case, the dominating term in \eqref{EqExpectDisc2} is the second one.
\end{rem}

The equality between the two lines of the theorem's equation \eqref{EqExpectDisc2} comes from the following elementary lemma.

\begin{lemma}\label{LemTheoDisc}
Let $\mu$ be a probability measure on $\Sp^1$, with distribution function minus average $f$. Let $F$ be a primitive of $f$ such that $F(0) = F(1) = 0$. Then
\[\int_0^1 (f^2+2F) = \Disc\big(\mu,\Leb)^2 - \frac{1}{12}\le 0.\]
\end{lemma}

A proof of this lemma can be found in appendix.

%
%

By taking $m_i=1$ for any $i$, one gets the following corollary about point precesses with points of uniform weights.

\begin{coro}\label{CoroDiscrLocDistrib}
Let $\mu$ be a probability measure on $\Sp^1$,  $N\in\N$, and $\nu_\mu$ be the random measure defined by
\[\nu_\mu = \frac{1}{N} \sum_{k=1}^N \delta_{p_k},\]
where the $p_k$'s are iid points with distribution $\mu$. Then 
\begin{equation}\label{EqExpectDisc}
\mathbb{E}\big[\Disc(\mu,\nu_\mu)^2\big] = \frac{1}{N} \left( \frac{1}{12} - \Disc\big(\mu,\Leb)^2\right).
\end{equation}
\end{coro}

\begin{rem}
A simple computation shows that the square of the distance $\Disc$ between $N$ equispaced points and Lebesgue measure, is equal to $1/(12N^2)$. On the other hand, one has $\mathbb{E}\big[\Disc(\Leb,\nu_{\Leb})^2\big] = 1/(12N)$. 

Hence, the expectation of the square of the Cramér distance of the uniform point process (with all weights equal to 1) is $1/(12N)$.
In this case, the squared distance $\Disc^2$ for a typical point process is way bigger than the minimal squared distance $\Disc^2$ for the same number of points (it is of the order of its square).
\end{rem}

\subsection{Comparison between the random point process and the discretization}

\begin{figure}
\begin{center}
\includegraphics[width=\linewidth]{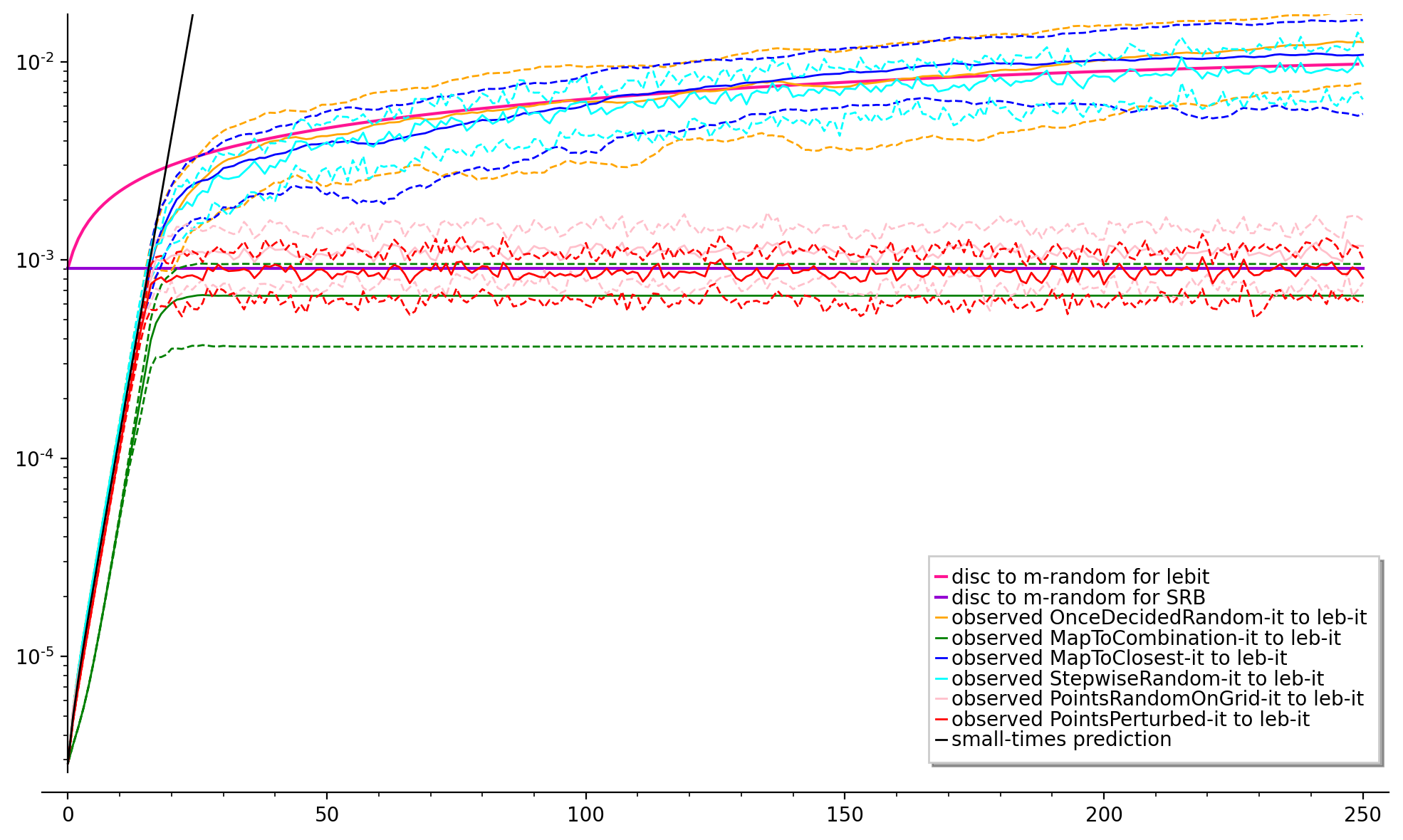}

\includegraphics[width=\linewidth]{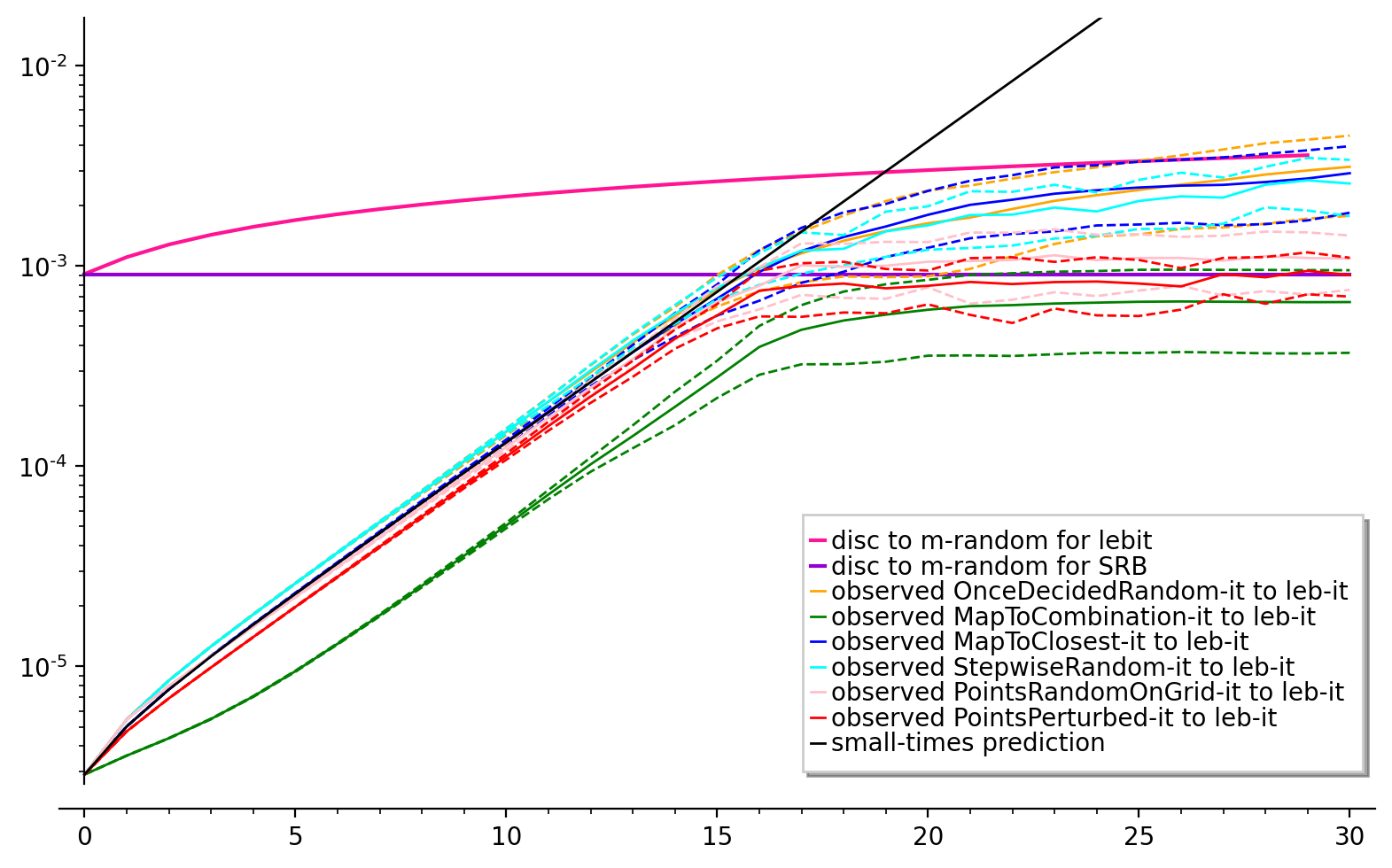}
\end{center}
\caption{\label{DiscrepancyIteratesDisc2} The bottom graphic is a zoom of the top one for small times. In these simulations, $N=10^5$. Each trio of curves of the same color (one in plain line and two dashed) is similar to the curves of Figure~\ref{FigSmallTimes}: it represents the mean (plain) and the mean $\pm$ standard deviation (dashed) of the distance $\Disc$ between the iterates of the discrete measure by the discretization type and the iterate of Lebesgue measure by RPF operator. The black curve is, as the red curve of Figure~\ref{FigSmallTimes}, the theoretical prediction given by Theorem~\ref{MainTheo}. The violet curve "disc to m-random for SRB" is the expected value of the distance between the SRB measure and $N$ random points with respect to this SRB measure, it is given by Corollary~\ref{CoroDiscrLocDistrib}. The pink curve "disc to m-random for lebit" is obtained from Theorem~\ref{TheoExpectDisc} (see the map $p(k)$ defined in Equation~\eqref{EqDefPk}).}
\end{figure}

Now we have defined different discretization types in Paragraph~\ref{ParagDifferent} and got a theoretical estimate of the distance $\Disc$ between a measure and the random point process associated to it paragraph~\ref{ParagPointProcess}, we can compare them numerically.

Figure~\ref{DiscrepancyIteratesDisc2} displays the mean values as well as the mean values $\pm$ standard deviation of the Cramér distance between the iterates of the discrete measure by the discretization type and the iterate of Lebesgue measure by RPF operator, for all the discretization types defined in Paragraph~\ref{ParagDifferent}.

It also shows the the theoretical prediction for small times given by Theorem~\ref{MainTheo}, and the curves of two expected values of the distance $\Disc$ between SRB measure and a point process. The first one, in violet, is obtained from Corollary~\ref{CoroDiscrLocDistrib} ; it is square root of the expected value of the square of the distance between the SRB measure and the measure made of $N$ independent random points with respect to this SRB measure. The second one, in pink, is cooked from Theorem~\ref{TheoExpectDisc} and the local distribution of preimages given by Proposition~\ref{PropLocDistrib1} in the following way: for each time $k$, we compute the theoretical local distribution of preimages $a_i=a_i(y,k)$ by the algorithm given by Proposition~\ref{LemLocDistrib}. To each of these functions $a_i(y,k)$, which represent local densities of measures (which we normalize to get probability measures), are associated a repartition functions $\tilde f_i^k$. This allows to get an estimation by means of \eqref{EqExpectDisc2} (see also Remark~\ref{RemExpectDisc})
\begin{equation}\label{EqDefPk}
p(k) = \sum_{i=1}^K \left(\int_{\Sp^1} a_i(y,k)\ud y\right)^2 \left(\frac{1}{12} - \Disc\left(\tilde f_i^k,\, \Leb\right)^2\right).
\end{equation}
\bigskip

At first glance, we can group the discretization types in three different clusters.
\begin{enumerate}
\item[(C1)] A first one containing only \texttt{MapToCombination}, whose asymptotic behaviour looks stationary, the asymptotic average distance is smaller than the estimation \texttt{disc to m-random for SRB}.
\item[(C2)] A second one containing \texttt{PointsRandomOnGrid} and \texttt{PointsPerturbed}. These two types of discretization give asymptotic behaviours similar to the estimation \texttt{disc to m-random for SRB}.
\item[(C3)] A third one containing \texttt{MapToClosest}, \texttt{OnceDecidedRandom} and \texttt{StepwiseRandom}. These three types of discretization behave asymptotically   more or less as \texttt{disc to m-random for lebit} given by the map $p(k)$.
\end{enumerate}

However, a closer look at each cluster reveals small differences.

For (C2), while the average value for \texttt{PointsPerturbed} follows asymptotically very well the curve of \texttt{disc to m-random for SRB}, the discretization type  \texttt{PointsRandomOnGrid} has a significantly greater average asymptotic value. This is quite unexpected, as the difference induced by replacing the SRB measure by the projection of it on $E_N$ (by mean of $E_N$) in Corollary~\ref{CoroDiscrLocDistrib} is of order $1/N^2$, and hence is negligible with respect to the orders of the computed distances, which are of order $1/N$. We have no explanation to this phenomenon.

\begin{figure}
\begin{center}
\includegraphics[width=\linewidth]{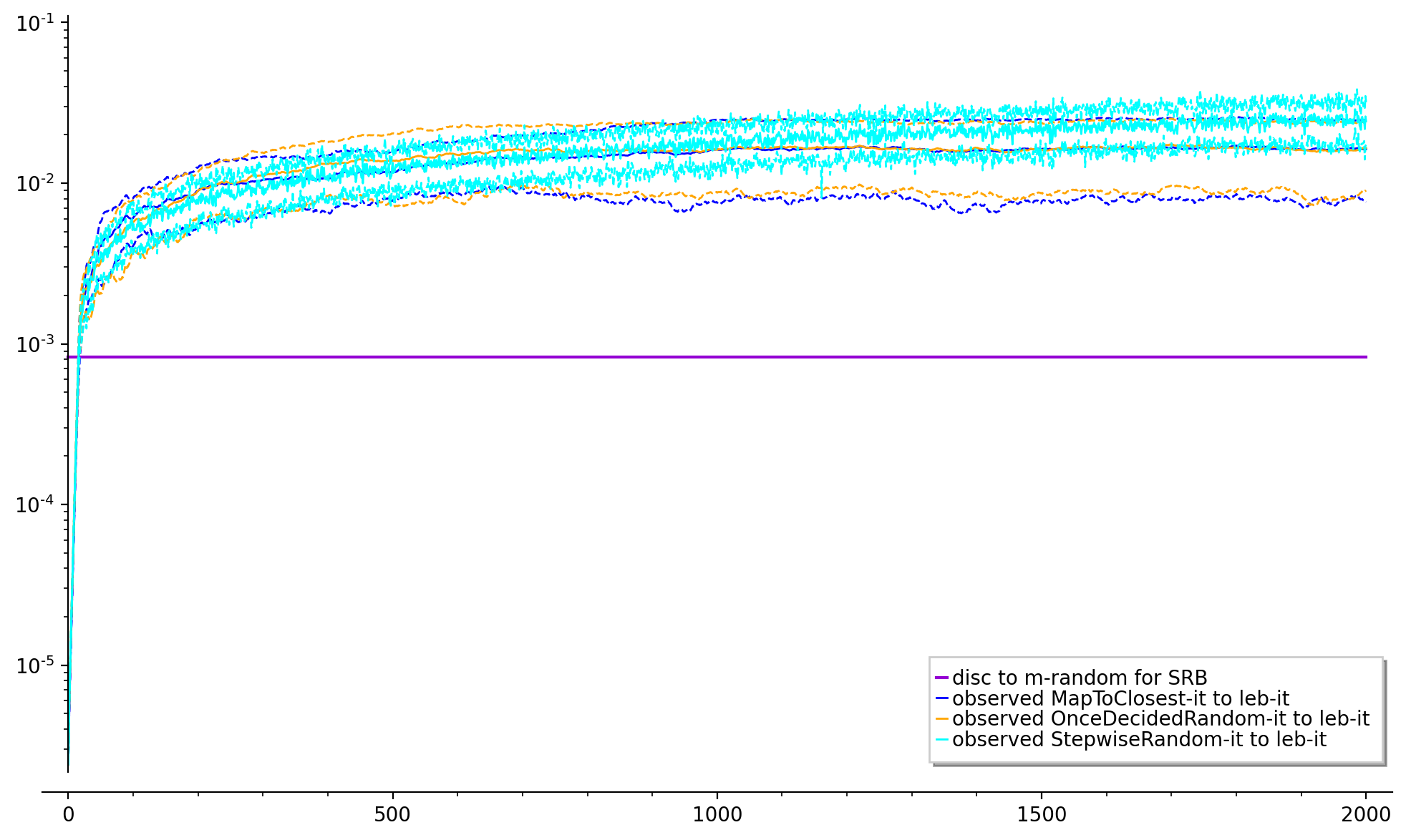}
\end{center}
\caption{\label{DiscrepancyIteratesDisc4} This figure shows the same curves as Figure~\ref{DiscrepancyIteratesDisc2} for $N = 120\,000$ but only for \texttt{MapToClosest}, \texttt{OnceDecidedRandom} and \texttt{StepwiseRandom}.}
\end{figure}

For (C3), the curves associated to \texttt{MapToClosest} and \texttt{OnceDecidedRandom} are very similar (see also Figure~\ref{DiscrepancyIteratesDisc4}). Hence, the already discussed difference between microscopic behaviours of these discretizations, the first one having deterministic local correlations and the second one random local correlations, seems to have no impact on the asymptotics of the distance $\Disc$ between measures. Moreover, the curves corresponding to \texttt{StepwiseRandom}, although also following quite well the curve of $p(k)$, behaves a bit more erratically than the two other ones (this is even more blatent on Figure~\ref{DiscrepancyIteratesDisc}).

Simulations over a larger time range show that the average distance for \texttt{StepwiseRandom}, form a certain point, gets bigger than the one for \texttt{MapToClosest} and \texttt{OnceDecidedRandom}. More precisely, Figure~\ref{DiscrepancyIteratesDisc4} shows this distance for times $t\le 2000$ and $N=120\,000$. On this simulation, we can see that this time where the mean distances start to be different is more or less $1000$. This time is to be compared with the mean time necessary for an orbit to cycle for a typical map of $120\,000$ elements, which is $\sqrt{\pi N/2}\simeq 137$: from this time 137, a significant part of orbits have cycled, which perturbs the process of injectivity loss. Note that on Figure~\ref{DiscrepancyIteratesDisc4} we have not represented the prediction that we discussed in the beginning of this paragraph, as in this time range accurate computations are out of our machine capacities: in the simulations we have to truncate the series $a_i$ up to some $i\le i_0$ to avoid exponential explosion of data depending on simulation time; we checked empirically that this truncation does not affect the prediction by verifying that the predictions are the same weather we truncate up to $i_0$ or $2i_0$. The threshold we chose for the simulations ($i_0=256$) gives similar results to $2i_0=512$ for times $\le 250$ (as in Figure~\ref{DiscrepancyIteratesDisc2}) but not $\le 2\,000$ (as in Figure~\ref{DiscrepancyIteratesDisc4}); a larger threshold for time $2\,000$ would make the computations extremely long and memory costing.
\bigskip

From these observations, one can conclude the following moral.

\begin{moral}
\textsf{The main phenomenon influencing the middle term behaviour of $(f_N^k)^*(\Leb_N)$ is the fact that orbits under $f_N$ merge.
More precisely, the distance $\Disc$ between $(f^k)^*(\Leb)$ and $(f_N^k)^*(\Leb_N)$ is rather well described by the distance $\Disc$ between the SRB measure and the point process described by $p(k)$ (see \eqref{EqDefPk}): locally around $y\in\Sp^1$, the proportion of points with weight $i$ of this process is equal to $a_i(y,k)$, where $a_i(y,k)$ represents the local proportion of points around $y$ that have $i$ preimages under $f_N^k$.}
\end{moral}

Note that in the simulations we performed the ``middle term'' is not that long and it may be that the phenomena specific to the short and long term still interfere in the time range and the discretization orders we chose: as noticed in Figure~\ref{DiscrepancyIteratesDisc}, left, for smaller orders $N$ there is even no middle term transitory behaviour. We would need more computing power to test the validity of the prediction \eqref{EqDefPk} on bigger orders $N$ (typically close to $10^7$).

\section{Long term behaviour}\label{SecAsym}

From Figure~\ref{FirstSimul}, one can wonder whether the quantity
\[\limsup_{k\to +\infty} \Disc\big((f_N)_*^k(\Leb_N),SRB\big),\]
which depends on $N$, tends to 0 when the discretization parameter $N$ goes to infinity or not. According to the simulations (see Figure~\ref{FigDiscSRBLimit}), it seems that yes, but unfortunately the proof of this result seems unreachable for now.

\begin{figure}[ht]
\centering \includegraphics[width=.5\linewidth]{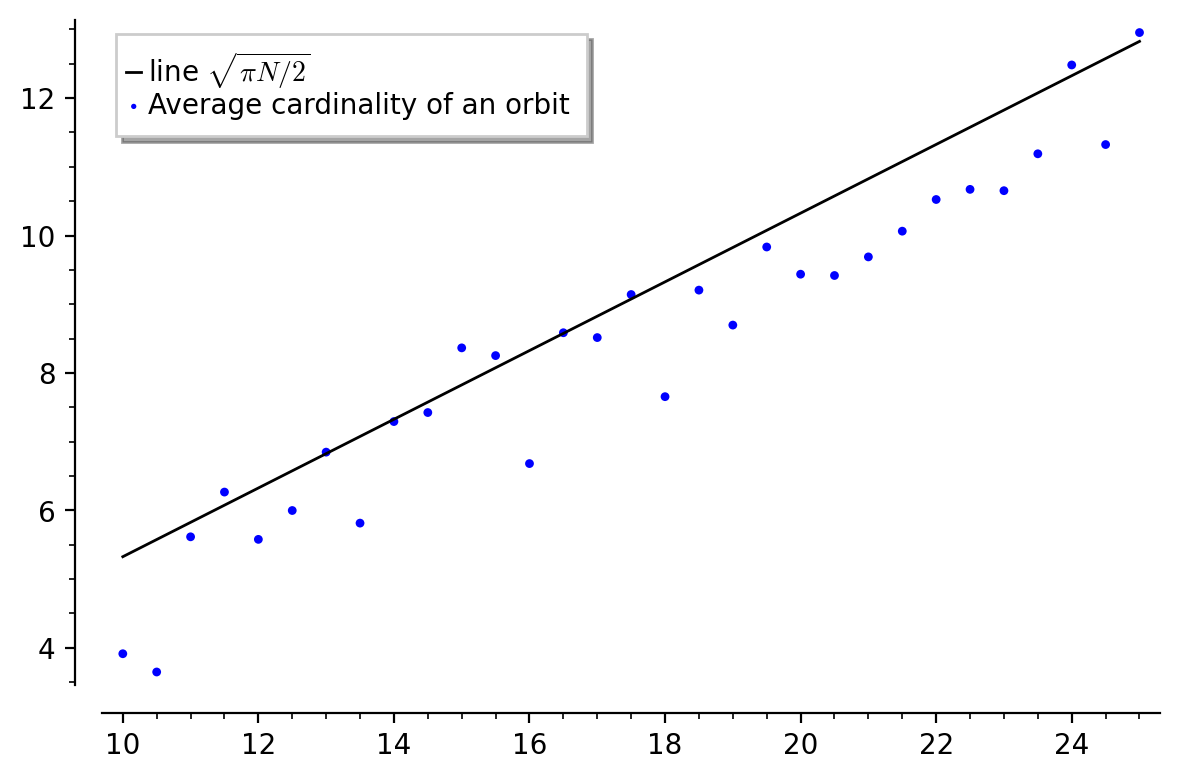}
\caption{\label{TimeToCycle} $\log_2$ of the average time needed for an orbit of $f_{2^k}$ to cycle (i.e. average cardinality of an orbit) depending on $k$ (blue points). This behaves quite the same as the same quantity for a typical random map on a set of $2^k$ elements (black line).}
\end{figure}

\begin{figure}[ht] 
\includegraphics[width=.49\linewidth]{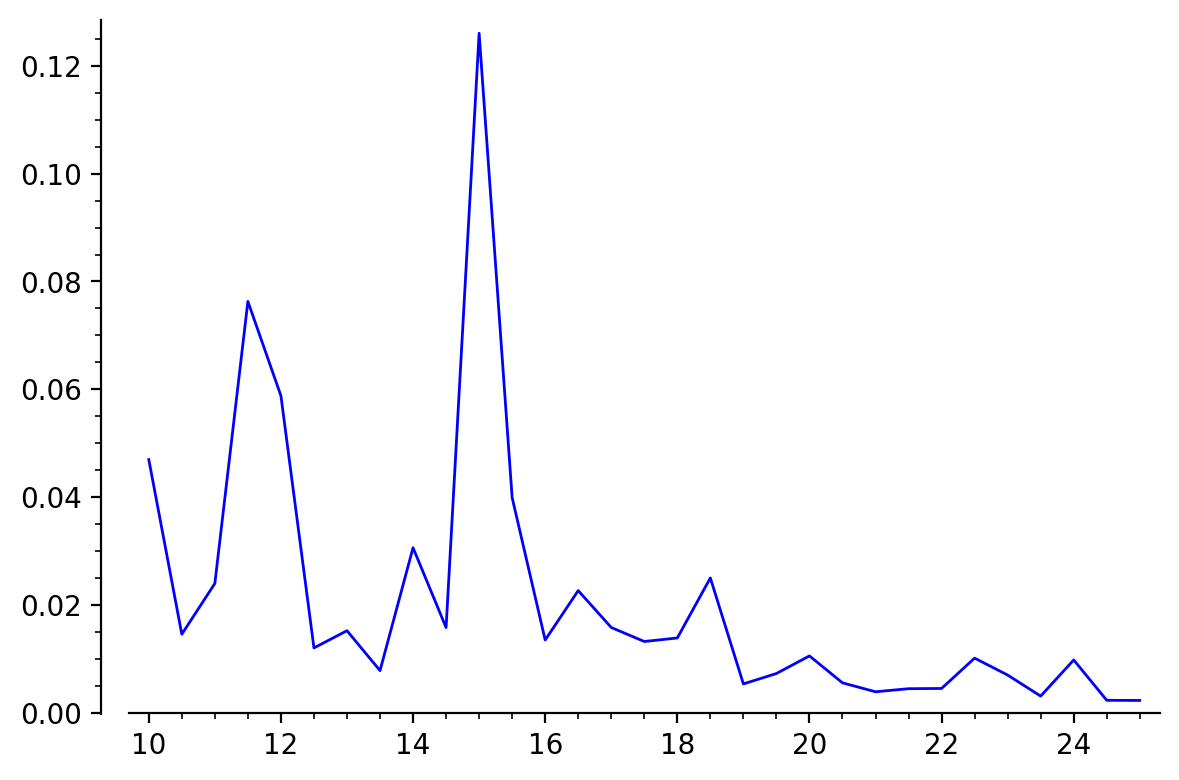}
\hfill
\includegraphics[width=.5\linewidth]{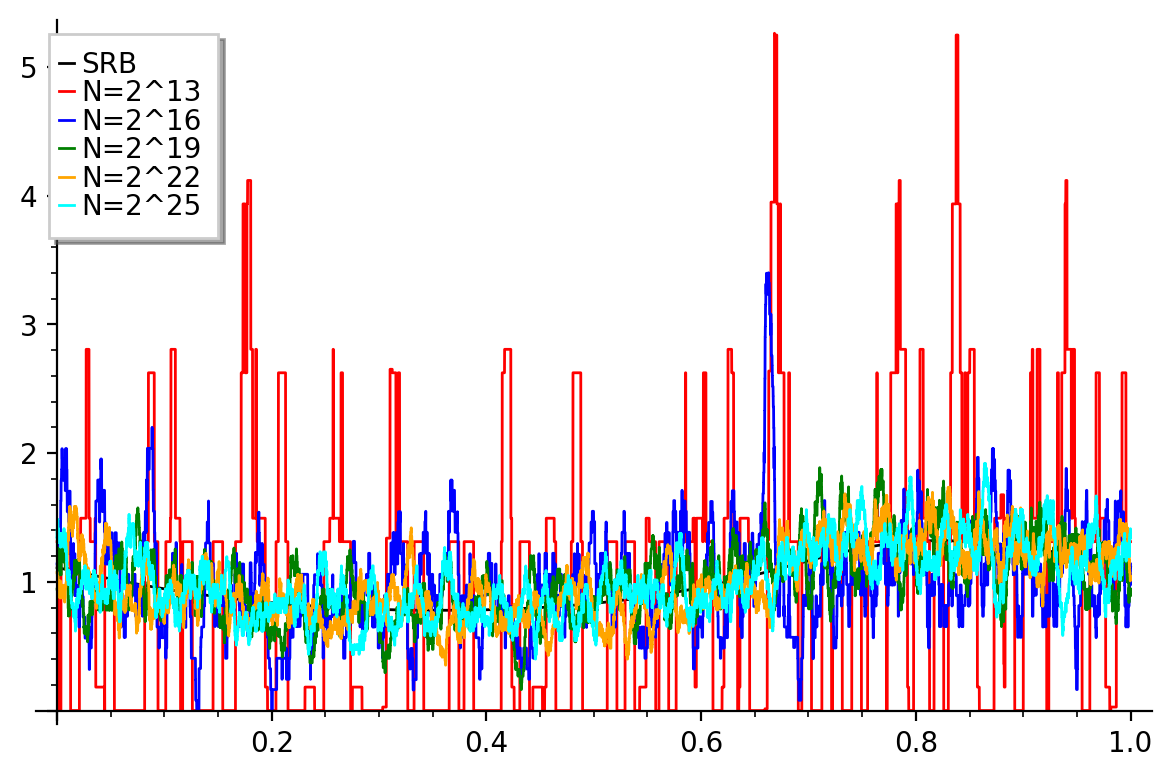}
\caption{\label{FigDiscSRBLimit} Left: Cramér distance between the measures $SRB$ and $\mu_{2^k}$ depending on $k$. We do not have any explanation for the threshold $k=15$. Right: densities of the measures $\mu_{2^k}$, for $k=13,16,19,22,25$, versus density of the SRB measure.}
\end{figure}

Let us first recall the fact that the maps $f_N$ are finite, so that every orbit eventually falls in a periodic cycle. Hence, one can characterize the expression ``asymptotic regime'' by the fact that all points have already browsed a whole periodic orbit\footnote{In practical, we will see the asymptotic regime's behaviour as soon as \emph{most of} points of $E_N$ already have browsed a whole cycle.}. Figure \ref{TimeToCycle} shows the average of this time over points of $E_{2^k}$ depending on $k$.

One can observe that this quantity behaves more or less as the same quantity for a typical random map on a set of $N$ elements, which is equivalent to $\sqrt{\pi N/2}$ (see \cite{Boll-rand}); this equivalent is represented in black in Figure~\ref{TimeToCycle}. This quantity is around $10^8$ for $N = 2^{56}$ (which is the classical precision used by computers). Hence, in practical, one usually does not reach the asymptotic regime when iterating a map.

There is a canonical measure associated to the asymptotic regime in the following way. Fix $N>0$, and let $\Leb_N$ be the uniform measure on the grid $E_N$. The measures $(f_N)_*^n(\Leb_N)$ converge in the Ces{\`a}ro mean towards a measure
\begin{equation}\label{EqDefMuN}
\mu_N = \lim_{k\to+\infty} \frac{1}{k}\sum_{n=0}^{k-1}(f_N)_*^n(\Leb_N).
\end{equation}
This measure is supported in the union of the periodic cycles of $f_N$, the total weight of each of them being proportional to the size of its basin of attraction.

\begin{ques}\label{ZeQuestion}
Do we have $\mu_N \to_{N\to +\infty} SRB$ (in the weak-* topology) for a generic $C^r$ expanding map $f$?
\end{ques}

Numerical experiments suggest that the answer to this question might be yes, as shown by Figure~\ref{FigDiscSRBLimit}. Note that the convergence, if happens, is extremely slow: it would give a terrible algorithm to compute an approximation of the SRB measure.

\begin{figure}[ht] 
\centering \includegraphics[width=.6\linewidth]{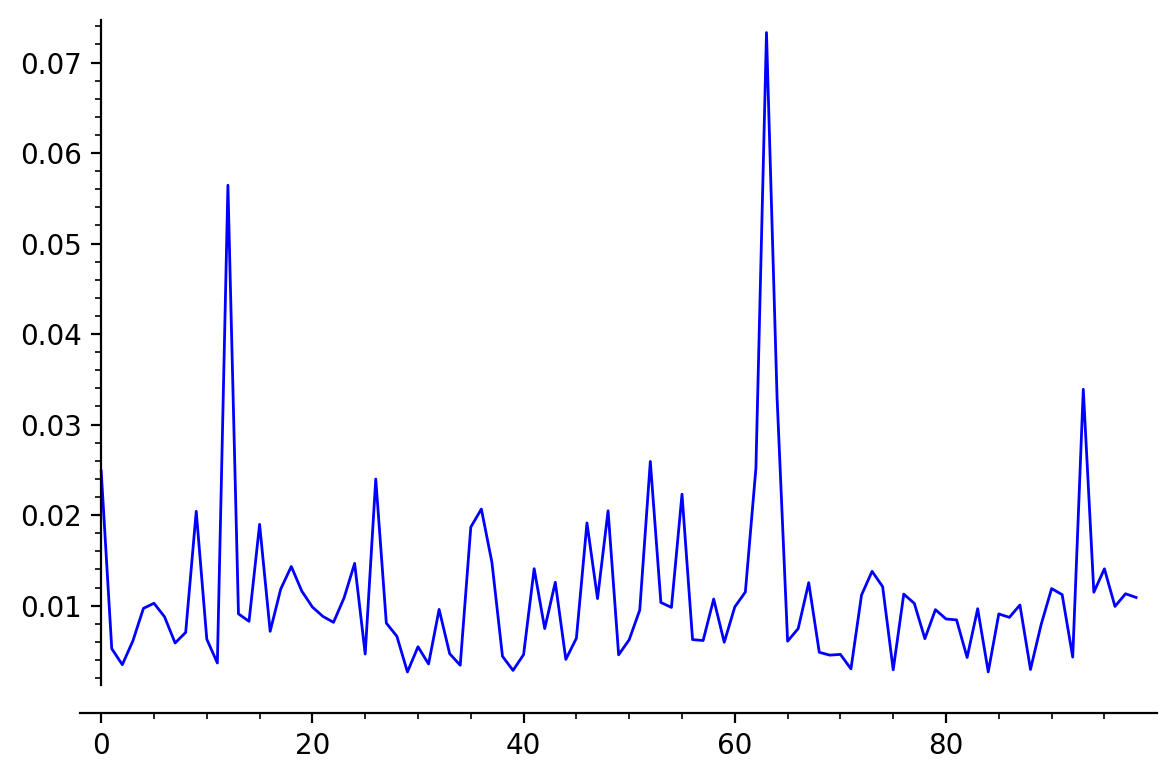}
\caption{\label{FigDiscSRBLimit2} Cramér distance between the measures $SRB$ and $\mu_{2^{20}+i}$, depending on $0\le i\le 99$.}
\end{figure}

However, it is possible that these simulations are misleading. As observed in \cite{Guih-These} (see Figures~12.14 and 12.17), for some $C^1$ area-preserving diffeomorphisms of the torus, the measures $\mu_N$ seem to converge to the area for a set of $N$ of density 1, but there are still rare values of $N$ for which the distance between $\mu_N$ and the area stay at positive distance. Such a behaviour can be observed on Figure~\ref{FigDiscSRBLimit2}: while for most of the  100 discretization orders $N$ between $2^{20}$ and $2^{20}+99$, the distance between the SRB measure and $\mu_N$ is around $0.01$, for two of these discretization orders, the distance is bigger than $0.05$ (a distance that is no longer attained after $2^{15}$ on Figure~\ref{FigDiscSRBLimit}).

\begin{moral}
\textsf{Our simulations do not suggest a clear conjecture about the convergence or not of $\mu_N$ towards $SRB$, but it seems that this convergence holds for a subsequence of $N$ of density 1.}
\end{moral}

\section{Conclusion}

We have identidfied three different temporal regimes and, for each one, proposed a model to describe it.

The first one, \ref{R1}, studied in Section~\ref{SecShort}, is the short term. It occurs for times $t\ll \log N$. In this regime, iterates of points that are initially microscopically close stay at a microscopical distance one from the others. The evolution of \eqref{EqStudiedQuant} is well described by Theorem~\ref{TauxExpand}.

The second one, \ref{R2}, is the middle term, studied in Section~\ref{SecMiddle}. It seems to occur for times $\log N\ll t$ and $\log t\ll \log N$. In this regime, orbits of points that were neighbours at time 0 are typically at a macroscopic distance one from the other.
Our simulations suggest that the main phenomena governing the evolution of \eqref{EqStudiedQuant} in this regime is \ref{P2}: two points of the grid $E_N$ having the same image by $f_N$ will have identical positive orbits.
On our simulations, the behaviour of \eqref{EqStudiedQuant} looks a lot like the one of a random process described in Paragraph~\ref{ParagPointProcess} (see \eqref{EqDefPk}). This is a point process made of points with different weights, the local density of the points with weight $p\in\N$ being equal to the predicted local density of points with $p$ preimages under the discretization $f_N^t$, this prediction being made by Proposition~\ref{PropLocDistrib1} (which reflects Phenomenon \ref{P2}). The expected Cramér distance $\Disc$ between the point process and the SRB measure is given by Theorem~\ref{TheoExpectDisc}.

The third one, \ref{R3}, the long term, is studied in Section~\ref{SecAsym}. It should happen for $\log t \gg \log N$ (or maybe $t\gg \sqrt N$, or $t\gg N$). What is the right model for this regime is more unclear than for the two other ones. It seems like that, as observed in the litterature (e.g. \cite{Flocker, MR1678095,MR701179,Mier-dyna, MR1400185, MR1169615, MR646380}), the combinatorial behaviour is well described by the one of a random map on a set with $N$ elements. Note that for some classes of interval maps having 0 as a fixed point, the combinatorial behaviour seem to be well described by a specific type of random mpas called \emph{random maps with a single attractive centre} \cite{MR1353178, MR1392078}. Even if our simulations suggest that the asymptotic measures $\mu_N$ may converge towards $SRB$ when $N$ goes to $+\infty$, this conclusion is not so clear in view of phenomena such as rare orders $N$ for which the measure $\mu_N$ is for away from SRB \cite{Guih-discr,Guih-These}.
\bigskip

Of course, it is natural to ask whether such numerical phenomena hold in higher dimensions for systems with some hyperbolic properties (Anosov or Axiom A systems, systems with dominated splittings, etc.).

It is indeed a real challenge to tackle a theoretical validation of the numerical observations we outlined in this paper. As written by Lanford in \cite{MR1678095}, ``\emph{[\dots] this problem may be as hard of that of non-equilibrium statistical mechanics.}''

We woulk like to put forward some research tracks that may be the first ones to address. First, one could try to get explicit times of convergence in Theorems~\ref{MainTheo} and \ref{TauxExpand} for some specific examples of piecewise real-analytic maps. Going a bit further, understanding why the first one is valid only in the short term and the second one remains true in the medium term would be an exceptionnal progress. 

Also, adapting the proofs of \cite{paper1} to the case of infinite branches circle expanding maps (e.g. Gauss map) may suggest other research directions. The program of previous paragraph may also be addressed in the case of the $\beta$-shift ($x\mapsto\beta x\mod 1$):  the constant slope of the map may allow to undestand completely the discretizations' behaviours.

%
%
%
%
%
%
%
%


\appendix

\small

\section{Proof of Theorem \ref{TheoExpectDisc} and Lemma \ref{LemTheoDisc}}

We start with the proof of Lemma \ref{LemTheoDisc}, some facts of it being used in the one of Theorem \ref{TheoExpectDisc}.

\begin{proof}[Proof of Lemma~\ref{LemTheoDisc}]
One has
\begin{align*}
\Disc\big(\mu,\Leb)^2
& =\int_0^1 \left(f(x)+\frac12-x\right)^2\ud x\\
& = \int_0^1 f^2 + 2\int_0^1\left(\frac12-x\right)f(x)\ud x +\int_0^1\left(\frac12-x\right)^2 \ud x\\
& = \int_0^1 f^2 - 2\int_0^1 xf(x)+\frac{1}{12}\\
& = \int_0^1 (f^2+2F) +\frac{1}{12}\qquad \text{(by parts),}
\end{align*}
so 
\[\int_0^1 (f^2+2F) = \Disc\big(\mu,\Leb)^2 - \frac{1}{12}.\]
This proves the equality of the lemma.

We now prove the inequality. We first suppose that $f$ is a convex combination of maps $g_p$, where $g_p(x) = \chi_{[p,1]}(x)-(1-p)$ is the cumulative-minus-average distribution function of $\delta_p$. For $g_p$ we have the equality
\[\langle g_p,g_p\rangle = p(1-p) = -2\int G_p.\]
If $f = \sum_i \lambda_i g_{p_i}$, with $p_i\in\Sp^1$ and $\lambda_i\ge 0$, $\sum_i\lambda_i = 1$, then
\begin{align*}
\int(f^2+2F) & = \big\langle \sum_i \lambda_i g_{p_i}, \, \sum_j \lambda_j g_{p_j}\big\rangle + 2\sum_i\lambda_i \int G_{p_i}\\
& = \sum_{i,j}\lambda_i\lambda_j\langle g_{p_i}, \, g_{p_j}\rangle - \sum_i\lambda_i \langle g_{p_i},\, g_{p_i}\rangle\\
& = \sum_{i,j}\lambda_i\lambda_j\big\langle g_{p_i}, \, g_{p_j}-g_{p_i}\big\rangle.
\end{align*}
Notice that fixing $p\le q$, we have
\begin{align}
\langle g_p,g_q \rangle &= \int_0^p (p-1)(q-1) \ud x + \int_p^q p(q-1) + \int_q^1 pq \nonumber \\
                            &= p(p-1)(q-1) + (q-p)p(q-1) + (1-q)pq \nonumber \\
                            &= p(1-q), \label{Eqgigj}
\end{align}
so for any $p,q$ we have $\langle g_p,g_q \rangle \ge \langle g_p,g_p \rangle$.
This proves that $\int(f^2+2F)\le 0$ in the case $f$ is a convex combination of maps $g_p$. The general case comes from the density (in $L^2$ norm) of such convex combinations among zero-average maps.

Note that a simple additional argument shows that we have equality $\int(f^2+2F) = 0$ iff $f=g_p$ for some $p\in\Sp^1$.
\end{proof}

\begin{proof}[Proof of Theorem \ref{TheoExpectDisc}]
Given two probability measures $\mu$ and $\nu$, the square $\Disc^2$ of their Cramér distance is obtained as
\[\int_0^1 (f(x)-g(x))^2 \ud x\]
where $f, g$ are the cumulative distribution function minus their average (as defined before Theorem \ref{TheoExpectDisc}).
In our case we want to compute the expectation of the square $\Disc^2$ of the Cramér distance between the measure $\mu$ and the random measure $\nu_\mu$ (which depends on the random points $p_i$). In this case, the (random) function $g$ is given by
\[g(x) = \frac{1}{N}\sum_{i=1}^{K}m_i g_{p_i}(x) = \sum_{i=1}^{K}\lambda_i g_{p_i}(x),\]
where $g_p(x) = \chi_{[p,1]}(x)-(1-p)$ is the cumulative-minus-average distribution function of $\delta_p$, and $\lambda_i = \frac{m_i}{N}$ (note that they satisfy $\sum_i \lambda_i = 1$). Each point $p_i$ will be chosen randomly and independently with distribution $\tilde \mu_{m_i}$
\bigskip

First, suppose that the distribution functions $\tilde f_m$ of $\tilde \mu_{m}$ are differentiable; in this case $\tilde f_m'$ is equal to the density of $\tilde \mu_{m}$. Note that in this case $\tilde f_m'\in L^1(\Sp^1)$, so that all further applications of Fubini's theorem will be valid.

Keeping implicit that we will be averaging over $p_1,\dots,p_n$ chosen at random in
$[0,1]^n$ with distribution $\tilde f_{m_1}'(p_1)\ud p_1 \cdots \tilde f_{m_n}'(p_n)\ud p_n $, the square $D$ of the Cramér distance $\Disc^2$ can be computed as
\begin{align*}
D = & \int_0^1 \left(f(x)-\sum_{i=1}^K \lambda_i g_{p_i}(x)\right)^2\ud x \\
 = &  \int_0^1 f(x)^2\ud x- 2\int_0^1\sum_{i=1}^K f(x)\lambda_i g_{p_i}(x) \ud x \\
 & +\int_0^1\sum_{i=1}^K \lambda_i^2 g_{p_i}(x)^2\ud x +\int_0^1\sum_{\substack{i,j=1\\i\neq j}}^K \lambda_i g_{p_i}(x)\lambda_j g_{p_j}(x)\ud x \\
= & \int_0^1 f(x)^2\ud x - 2\sum_{i=1}^K \lambda_i \int_0^1 f(x) g_{p_i}(x)\ud x \\
 & + \sum_{i=1}^K \lambda_i^2 \int_0^1 g_{p_i}(x)^2\ud x + \sum_{\substack{i,j=1\\i\neq j}}^K  \lambda_i \lambda_j \int_0^1 g_{p_i}(x) g_{p_j}(x)\ud x.
\end{align*}
In each integral the average will only be over $p_i$ (or $p_i, p_j$ in the last one), as the integrand does not depend on $p_k$, for $k\neq i,j$. Each $p_i$ follows the distribution $\tilde f_{m_i}'(p)\ud p$. Therefore we can simplify the sums of equal values and see the integrals as integral in just $p$ (or $p,q$ in the last integral). We can write:
\begin{align}
D = & \int_0^1 f(x)^2\ud x - 2\sum_i \lambda_i\int_{p=0}^1 \int_{x=0}^1 f(x)g_{p}(x) \tilde f_{m_i}'(p)\ud x\ud p \nonumber \\
  & + \sum_i \lambda_i^2 \int_{p=0}^1\int_{x=0}^1 g_{p}(x)^2 \tilde f_{m_i}'(p)\ud x\ud p \nonumber\\
  & + \sum_{i\neq j}\lambda_i\lambda_j \int_{p=0}^1\int_{q=0}^1 \int_{x=0}^1 g_{p}(x)g_{q}(x)\tilde f_{m_i}'(p) \tilde f_{m_j}'(p)\ud x \ud p \ud q\nonumber\\
= & I_1 - 2 I_2 + \sum_i \lambda_i^2\, I_3 + \sum_{i\neq j} \lambda_i\lambda_j I_4 \label{eq1}
\end{align}

Before deducing a general formula we well state two trivial lemmas.

\begin{lemma} \label{flemma1}
For all $a\leq b$ we have:
  \begin{align*}
  \int_a^b pf'(p)\ud p &= \big[pf(p) \big]_a^b - \int_a^b f(p)\ud p \\
                       &= bf(b)-af(a) - F(b) + F(a).
\end{align*}
In particular, we get $\int_0^b pf'(p)=bf(b)-F(b)$.
\end{lemma}

\begin{lemma} \label{flemma2}
\begin{align*}
\int_a^b f'(p)F(p)\ud p &= \big[f(p)F(p) \big]_a^b - \int_a^b f(p)^2\ud p.
\end{align*}
In particular we get $\int_0^1f'(p)F(p)=-\int_0^1 f(p)^2$.
\end{lemma}

The expression \eqref{Eqgigj} is valid on the half of the square $(p,q)\in[0,1]^2$ where $p\leq q$, therefore we have (applying the above lemmas and the fact that $F(0)=F(1)=0$)
\begin{align*}
I_4 = & \int_{p=0}^1 \int_{q=0}^1 \int_{x=0}^1 g_{p}(x)g_{q}(x)\tilde f_{m_i}'(p) \tilde f_{m_j}'(q)\ud x \ud q \ud p\\
= & \int_{q=0}^1\int_{p=0}^q p(1-q) \tilde f_{m_i}'(p)\tilde f_{m_j}'(q)\ud p \ud q + \int_{q=0}^1\int_{p=q}^1 q(1-p) \tilde f_{m_i}'(p)\tilde f_{m_j}'(q)\ud p \ud q \quad \text{(by \eqref{Eqgigj})} \\
= & - \int_{q=0}^1\int_{p=0}^1 pq \tilde f_{m_i}'(p)\tilde f_{m_j}'(q)\ud p \ud q\\
 & + \int_{q=0}^1\int_{p=0}^q p \tilde f_{m_i}'(p)\tilde f_{m_j}'(q)\ud p \ud q
 + \int_{p=0}^1\int_{q=0}^p q \tilde f_{m_i}'(p)\tilde f_{m_j}'(q)\ud q \ud p \\
= & - \left(\int_{p=0}^1 p \tilde f_{m_i}'(p)\right)\left(\int_{q=0}^1 q \tilde f_{m_j}'(q)\right)\\
 & + \int_{q=0}^1 \left(q \tilde f_{m_i}(q) - \tilde F_{m_i}(q)\right)\tilde f_{m_j}'(q) \ud q
 + \int_{p=0}^1 \left(p\tilde f_{m_j}(p) - \tilde F_{m_j}(p)\right)\tilde f_{m_i}'(p)\ud p \\
= & - \tilde f_{m_i}(1)\tilde f_{m_j}(1)\\
 & + \int_{q=0}^1 q\left(\tilde f_{m_i}(q)\tilde f_{m_j}'(q) + \tilde f_{m_j}(q)\tilde f_{m_i}'(q) \right) \ud q - \int_{q=0}^1 \left(\tilde F_{m_i}(q)\tilde f_{m_j}'(q) + \tilde F_{m_j}(q)\tilde f_{m_i}'(q)\right) \ud q \\
= & - \tilde f_{m_i}(1)\tilde f_{m_j}(1) \qquad\qquad \text{(by parts, for both integrals)}\\
 & + \left[ q \tilde f_{m_i}(q)\tilde f_{m_j}(q)\right]_0^1 - \int_{0}^1 \tilde f_{m_i}\tilde f_{m_j} - \left[ \tilde F_{m_i}\tilde f_{m_j} + \tilde F_{m_j}\tilde f_{m_i}\right]_0^1 + \int_{0}^1 \left(\tilde f_{m_i}\tilde f_{m_j} + \tilde f_{m_j}\tilde f_{m_i}\right)\\
= & \int_{0}^1 \tilde f_{m_i}\tilde f_{m_j} \numberthis \label{eqI4}
\end{align*}
For $I_3$ of \eqref{eq1} we have
\begin{align*}
I_3 & = \int_{p=0}^1\int_{x=0}^1 g_{p}(x)^2 \tilde f_{m_i}'(p)\ud x\ud p\\
 &= \int_{p=0}^1 p(1-p) \tilde f_{m_i}'(p)\ud p \\
   &= \big[p(1-p)\tilde f_{m_i}(p)\big]_0^1 + \int_0^1 2p \tilde f_{m_i}(p)\ud p - \int_0^1 \tilde f_{m_i}(p)\ud p\quad \text{(by parts)}\\
   &= \big[2p\tilde F_{m_i}(p)\big]_0^1 - \int_0^1 2\tilde F_{m_i}(p) \ud p\quad \text{(by parts)} \\
   &= - \int_0^1 2\tilde F_{m_i},\numberthis \label{eqI3}
\end{align*}
while for the second integral 
\begin{align*}
I_2 & = \int_{p=0}^1 \int_{x=0}^1 f(x)g_{p}(x) \tilde f_{m_i}'(p)\ud x\ud p\\
 & = \int_{p=0}^1 \int_{x=0}^1 f(x)\chi_{[p,1]}(x)\tilde f_{m_i}'(p)\ud x \ud p - \int_{p=0}^1 \int_{x=0}^1 f(x)(1-p)\tilde f_{m_i}'(p)\ud x \ud p \\
   & = \int_{x=0}^1 f(x)\int_{p=0}^x \tilde f_{m_i}'(p)\ud p \ud x - \int_{p=0}^1 \left( \int_0^1 f(x)\ud x \right) (1-p)\tilde f_{m_i}'(p)\ud p\\
   & = \int_{x=0}^1 \big(f(x)\tilde f_{m_i}(x) - f(x)\tilde f_{m_i}(0)\big) \ud x - 0 \quad \text{(because $f$ is zero-average)} \\
   & = \int_0^1 f\tilde f_{m_i} \numberthis \label{eqI2}
 \end{align*}
Joining all simplified expressions \eqref{eqI2}, \eqref{eqI3} and \eqref{eqI4}, we deduce that \eqref{eq1} can be rewritten as
\begin{align*}
D = & \int_0^1 f^2 - 2\sum_i \frac{m_i}{n} \int_0^1 f \tilde f_{m_i} - 2\sum_i \left(\frac{m_i}{n}\right)^2 \int_0^1\tilde F_{m_i} + \sum_{i\neq j} \frac{m_i m_j}{n^2} \int_0^1\tilde f_{m_i}\tilde f_{m_j}\\
 = & \int_0^1\left(f-\sum_i \frac{m_i}{n}\tilde f_{m_i} \right)^2 - \sum_i \frac{m_i^2}{n^2} \int\left(\tilde f_{m_i}^2 + 2\tilde F_{m_i}\right),
\end{align*}
which gives the first formula of the theorem. The second one is a consequence of Lemma~\ref{LemTheoDisc}
\bigskip

We now treat the general case for $\tilde f_m$. The cadlag map $\tilde f_m$ can be approached in uniform topology by a smooth repartition function of a measure $\overline\mu_m$, which is close to $\tilde \mu_m$ in weak-* topology. It then suffices to remark that in \eqref{EqExpectDisc}, the left side is continuous in $\tilde \mu_m$ for weak-* topology, and right side is continuous in $\tilde f_m$ for the uniform topology.
\end{proof}

\bibliographystyle{amsalpha}
\bibliography{../../Biblio.bib}

\vfill

\end{document}